\documentclass[amsfonts,11pt]{amsart}
\RequirePackage{amsmath,amssymb,amsthm, amscd, comment,mathtools}
\usepackage[margin=1.2in]{geometry}
\usepackage{amsmath}
\usepackage{amsthm}
\usepackage{amssymb}
\usepackage{amscd}
\usepackage{mathrsfs}
\usepackage{graphicx}
\usepackage{rotating, graphicx}
\usepackage[all,cmtip]{xy}
\usepackage{enumerate}
\usepackage{tikz-cd}
\usepackage[thinlines]{easytable}
\usepackage{multirow}
\usepackage{mathtools}
\usepackage{mathrsfs}
\usepackage{changepage} 
\usepackage{paralist}

\usepackage{kantlipsum}
\usepackage{enumitem}
\usepackage{times}
\usepackage{subfig}
\usepackage{color}
\usepackage[pagebackref,breaklinks,colorlinks,linkcolor=red,anchorcolor=red,citecolor=blue]{hyperref}

\calclayout

\newtheorem{proposition}[subsubsection]{Proposition}
\newtheorem{corollary}[subsubsection]{Corollary}
\newtheorem{theorem}[subsubsection]{Theorem}
\newtheorem{lemma}[subsubsection]{Lemma}
\theoremstyle{definition}
\newtheorem{definition}[subsubsection]{Definition}
\theoremstyle{remark}

\newtheorem{remark}[subsubsection]{Remark}
\newtheorem{example}[subsubsection]{Example}

\numberwithin{equation}{subsubsection}

\DeclareMathAlphabet{\mathbbold}{U}{bbold}{m}{n}

\title{Trace maps in motivic homotopy and local terms}

\author{Fangzhou Jin}
\address{School of Mathematical Sciences\\
Tongji University\\
Siping Road 1239\\
200092 Shanghai\\
P.R.China}
\email{\href{mailto:fangzhoujin@tongji.edu.cn}{fangzhoujin@tongji.edu.cn}}
\urladdr{\url{https://fangzhoujin.github.io/}}

\date{\number\day-\number\month-\number\year}

\begin{document}

\maketitle

\begin{abstract}

We define a trace map for every cohomological correspondence in the motivic stable homotopy category over a general base scheme, which takes values in the twisted bivariant groups. Local contributions to the trace map give rise to quadratic refinements of the classical local terms, and some $\mathbb{A}^1$-enumerative invariants, such as the local $\mathbb{A}^1$-Brouwer degree and the Euler class with support, can be interpreted as local terms. We prove an analogue of a theorem of Varshavsky, which states that for a contracting correspondence, the local terms agree with the naive local terms.

\end{abstract}

\tableofcontents

\noindent

\section{Introduction}
In Grothendieck's approach to the Weil conjectures via the theory of \'etale cohomology, a \emph{trace formula}, analogous to the Lefschetz fixed-point theorem in algebraic topology, plays a crucial role: it states that the number of rational points of an algebraic variety over a finite field can be computed as the alternating sum of the traces of the Frobenius morphism acting on (compactly supported) \'etale cohomology groups (\cite{Gro}). This formula is later generalized to what is now called the \emph{Lefschetz-Verdier formula} (\cite[III Cor. 4.7]{SGA5}), which states the proper covariance of a cohomological pairing called the \emph{Verdier pairing}. In both cases, the computation of the local contributions to the global trace is an interesting but difficult problem. The case of curves has been discussed in \cite[IIIb]{SGA5}; in the topological case, the work of \cite{GM} provides a practical formula for a class of maps called \emph{weakly hyperbolic}; over finite fields, the problem is studied at length in \cite{Pin} and \cite{Fuj}, where a conjecture of Deligne is eventually proved; in \cite{Var} the computation of local terms is generalized to \emph{contracting correspondences}.

The \emph{trace} of a cohomological correspondence over a field in the motivic context has been considered in \cite{Ols} and \cite{Cis} for \'etale motives, and generalized in \cite{JY} to the motivic stable homotopy category, assuming resolution of singularities. The trace map over a general base scheme was first introduced in \cite{YZ} and \cite{JY} under some smoothness and tranversality assumptions, and generalized in \cite{LZ} to the even singular base schemes.

Our first goal is to generalize such a construction over a field to a general base scheme:
\begin{theorem}[see \S~\ref{sec:corrtr}]
Let $C\xrightarrow{(c_1,c_2)}X\times_SX$ be a morphism of schemes, let $K\in\mathcal{SH}_c(X)$ be a constructible motivic spectrum over $X$ which is \emph{universally strongly locally acyclic} over $S$ (see Definition~\ref{def:locacy}), and let $v$ be a virtual vector bundle on $C$. Then for a \emph{cohomological correspondence} $u:c_1^*K\to c_2^!K\otimes Th(v)$ on can define its \emph{trace map}
\begin{align}
\label{eq:reltraceintro}
Tr(u/S):Th(-v_{|Fix(c)})\to \mathcal{K}_{Fix(c)/S}
\end{align}
where $Fix(c)$ is given by the fiber product
\begin{align}
\begin{split}
  \xymatrix@=10pt{
    Fix(c) \ar[r]^-{c'} \ar[d]_-{\delta'} & X \ar[d]^-{\delta_{X/S}}\\
    C \ar[r]^-{c} & X\times_SX.
  }
\end{split}
\end{align}
\end{theorem}
We use the techniques developed in \cite[\S2]{JY} and follow the ideas in \cite[\S2.1]{LZ}, by proving the K\"unneth formulas over a base in~\ref{sec:relkun}. In \S~\ref{sec:opcorr} we study the compatibility of the trace map with several operations: base change, push-forward and pullback.

Recent works on enumerative invariants given by intersection numbers in the framework of motivic homotopy theory (for example \cite{Hoy}, \cite{KW}, \cite{Lev}, \cite{LR}, \cite{BW}) have lead to a theory of \emph{$\mathbb{A}^1$-enumerative geometry}, where one define and study arithmetic refinements of classical invariants in terms of symmetric bilinear forms (or quadratic forms in characteristic different from $2$). 

Using the trace maps, in Section~\ref{sec:A1inv} we define the \emph{local terms} of the trace maps in motivic homotopy, which are local contributions to the trace map, and give quadratic refinements of the classical local terms. We show that some $\mathbb{A}^1$-enumerative invariants in the literature, such as the \emph{local $\mathbb{A}^1$-Brouwer degree} and the \emph{Euler class with support}, can be interpreted as local terms.

In the remaining of the paper, we prove the following analogue of Varshavsky's theorem (\cite[Thm. 2.1.3]{Var}):
\begin{theorem}[see Theorem~\ref{th:LTnaive}]
\label{th:intro}
Assume that the base $S$ is the spectrum of a field, and that the conditions~\ref{num:relori} and~\ref{num:Kcondfin} below are satisfied.
If the morphism $C\to X\times_SX$ \emph{contracting} (see Definition~\ref{def:invind} and Definition~\ref{def:contrnbh}), then the local terms agree with the naive local terms.
\end{theorem}
The \emph{naive local term} may be understood as follows, in the case where the correspondence is given by an endomorphism: if $X$ is a scheme, $f:X\to X$ is an endomorphism of $X$, $K\in\mathcal{SH}_c(X)$ is a constructible motivic spectrum over $X$, and $u:f^*K\to K$ is a map, then for every fixed point $x$ of $X$, the local term $LT_x(u)$, if it is well defined, is the trace of the induced map $u_x:K_{|x}\to K_{|x}$. 

The conditions~\ref{num:relori} are related to the problem of defining push-forward maps in motivic homotopy theory, which is related to the existence of \emph{orientations} and \emph{Thom classes}. See for example \cite[\S7]{DFJK}, \cite[\S4]{BW}. The conditions~\ref{num:Kcondfin} is related to the motivic nearby cycle functor, which will be used in the proof. They are satisfied for example if $k$ has characteristic $0$, or if we work with rational coefficients in $\mathcal{SH}_c(-,\mathbb{Q})$, see Lemma~\ref{lm:condsat} below.

The basic ideas of the proof are very close to \cite{Var}, with the following ingredients:
\begin{enumerate}
\item The use of \emph{additivity of traces} (\cite[Prop. 1.5.10]{Var}), which states that the trace map is additive along distinguished triangles. Note that the general statement fails for symmetric monoidal triangulated categories (\cite{Fer}), and the proof of \cite[Prop. 1.5.10]{Var} uses a variant for the filtered derived category. Our treatment here (Lemma~\ref{lm:additivity}) uses the language of higher categories, for which we refer to \cite[\S4]{JY} for a more detailed discussion.
\item The geometric fact that the under the deformation to the normal cone, a contracting correspondence becomes supported at the zero section of the normal cone (\cite[Rmk. 2.1.2 (b)]{Var}).
\end{enumerate}
The proof uses Ayoub's theory of \emph{motivic nearby cycles} (\cite[\S3]{Ayo}) which, together with the trace map over a base~\eqref{eq:reltraceintro}, allows to define the specialization of correspondences (\ref{sec:spcorr}). In Lemma~\ref{lm:speq} we show that this specialization is compatible with the \emph{specialization map} constructed in \cite[4.5.6]{DJK}, which is modeled on Fulton's specialization map on Chow groups (\cite[\S20.3]{Ful}). In Section~\ref{sec:A1Verdier} we prove an analogue of Verdier's theorem on restriction to vertices, where the proof is simplified compared to the classical case via the use of $\mathbb{A}^1$-homotopies and the description of generators of $\mathcal{SH}$ (Lemma~\ref{lm:CD7.2}). In Section~\ref{sec:defcorr}, we finish the proof of Theorem~\ref{th:intro} by putting up everything together.

\subsubsection*{Acknowledgments}
The author would like to thank Tom Bachmann, Fr\'ed\'eric D\'eglise, Adeel Khan, Marc Levine, Kirsten Wickelgren and Enlin Yang for helpful discussions. He would like to thank the referee for carefully reading the paper and giving suggestions that help improve the quality of the paper.
He acknowledges support from the National Key Research and Development Program of China Grant Nr.2021YFA1001400, the National Natural Science Foundation of China Grant Nr.12101455, the Fundamental Research Funds for the Central Universities, and the ERC Project-QUADAG, which has received funding from the European Research Council (ERC) under the European Union’s Horizon 2020 research and innovation programme (grant agreement Nr.832833). 

\section{Preliminaries}

\subsection{Notations and conventions}
\subsubsection{}
All schemes are assumed quasi-compact and quasi-separated, and all morphisms of schemes are separated of finite type. 
Regular schemes are assumed noetherian. 

\subsubsection{}
Throughout the paper, we denote by $\mathfrak{P}$ a set of prime numbers.
We will be mainly working with the motivic stable homotopy category $\mathcal{SH}(-)[\mathfrak{P}^{-1}]$, 
considered as a \textbf{motivic $\infty$-category} (\cite{Kha}) which is \textbf{continuous} (\cite[Appendix A]{DFJK}), and we denote by $\mathcal{SH}_c$ the full subcategory of \textbf{constructible objects} (\cite[Def. 4.2.1]{CD}). 

By \cite[Thm. 7.14]{DG}, it turns out that $\mathcal{SH}$ is the universal $\infty$-category which satisfies the six functor formalism, and it is possible to transport the main results to other motivic $\infty$-categories, which we do not discuss here.

\subsubsection{}
For a scheme $X$, we denote by $\mathbbold{1}_X$ the unit object of the symmetric monoidal $\infty$-category $\mathcal{SH}(X)[\mathfrak{P}^{-1}]$. For an object $K\in\mathcal{SH}(X)[\mathfrak{P}^{-1}]$ and an integer $n$, we denote by $K[n]=K\otimes (S^1)^{\otimes n}$ its $n$-th simplicial suspension.

\subsubsection{}
For $X$ a scheme, a \emph{virtual vector bundle} on $X$ corresponds to a point of the $K$-theory space of $X$. 
If $v$ is a virtual vector bundle on $X$, we denote by $Th(v)$ the associated \textbf{Thom space}, which is a $\otimes$-invertible object in $\mathcal{SH}(X)$ (\cite[\S16.2]{BH}). In particular, we define the \textbf{Tate twists} as $\mathbbold{1}_X(n)=Th(\mathbb{A}^n_X)[-2n]$.

\subsubsection{}
If $f:X\to Y$ is a morphism and $K\in\mathcal{SH}(Y)$, we denote $K_{|X}=f^*K\in\mathcal{SH}(X)$.

\subsubsection{}
For $f:X\to S$, let $\mathcal{K}_{X/S}=f^!\mathbbold{1}_S$, and 
\begin{align}
\begin{split}
\mathbb{D}_{X/S}:\mathcal{SH}(X)&\to\mathcal{SH}(X)\\
M&\mapsto \underline{Hom}(M,\mathcal{K}_{X/S}).
\end{split}
\end{align}

\subsubsection{}
\label{num:RS}
We say that a perfect field $k$ \textbf{satisfies strong resolution of singularities} if the following conditions hold:
    \begin{enumerate}
    \item For every separated integral scheme $X$ of finite type over $k$, there exists a proper birational surjective morphism $X'\to X$ with $X'$ regular;
    \item For every separated integral regular scheme $X$ of finite type over $k$ and every nowhere dense closed subscheme $Z$ of $X$, there exists a proper birational surjective morphism $b:X'\to X$ such that $X'$ is regular, $b$ induces an isomorphism $b^{-1}(X-Z)\simeq X-Z$, and $b^{-1}(Z)$ is a strict normal crossing divisor in $X'$.
    \end{enumerate}

\subsection{Bivariant theory in motivic homotopy}
\subsubsection{}
If $f:X\to S$ is a separated morphism of finite type and $v$ is a virtual vector bundle on $X$, define the \textbf{$v$-twisted bivariant spectrum} as the mapping spectrum
\begin{align}
H(X/S,v)=\operatorname{Map}_{\mathcal{SH}(X)}(Th(v),f^!\mathbbold{1}_S).
\end{align}
More generally, to every motivic spectrum $\mathbb{E}\in\mathcal{SH}(S)$ we associate a spectrum
\begin{align}
\label{eq:bivsp}
\mathbb{E}(X/S,v)=\operatorname{Map}_{\mathcal{SH}(X)}(Th(v),f^!\mathbb{E}).
\end{align}
If $\mathbb{E}$ is equipped with a unit section $\mathbbold{1}_S\to\mathbb{E}$, then it induces a map
\begin{align}
\label{eq:bivsec}
H(X/S,v)
\to
\mathbb{E}(X/S,v).
\end{align}

In what follows we recall its functorialities established in \cite{DJK}.
\subsubsection{}
\label{num:lci}
A morphism of schemes is said to be a \textbf{local complete intersection} if it can be (globally) factored as a regular closed immersion followed by a smooth morphism (\cite[\S6.6]{Ful}). If $f:X\to Y$ is a local complete intersection morphism, we denote by $\tau_f$ or $\tau_{X/Y}$ its virtual tangent bundle (\cite[Example 2.3.10]{DJK}), considered as a virtual vector bundle over $X$. Recall that one of the main results of \cite{DJK} states that there exists a natural transformation of functors
\begin{align}
\label{eq:purtr}
f^*(-)\otimes Th(\tau_f)\to f^!(-)
\end{align}
called \textbf{purity transformation}, which is induced by the \textbf{fundamental class} of $f$ 
\begin{align}
\eta_f\in H(X/Y,\tau_{X/Y})
\end{align}
(\cite[4.3.1]{DJK}). We say that an object $K\in\mathcal{SH}(X)$ is \textbf{$f$-pure} if the purity transformation $f^*K\otimes Th(\tau_f)\to f^!K$ is an isomorphism (\cite[Def. 4.3.7]{DJK}).

\subsubsection{Base change}(\cite[2.2.7 (1)]{DJK})
\label{num:DJK2271}
Consider a Cartesian square
\begin{align}
\label{eq:Cartsqdelta}
\begin{split}
  \xymatrix@=10pt{
    Y \ar[r]^-{q} \ar[d]_-{g} \ar@{}[rd]|{\Delta} & X \ar[d]^-{f}\\
    Z \ar[r]_-{p} & S
  }
\end{split}
\end{align}
and $v$ (respectively $w$) a virtual vector bundle on $X$ (respectively $Z$).
There are canonical maps 
\begin{align}
\label{eq:verBC}
\Delta^*:H(Z/S,w)\to H(Y/X,w_{|Y})
\end{align}
\begin{align}
\label{eq:horBC}
\Delta'^*:H(X/S,v)\to H(Y/Z,v_{|Y}).
\end{align}

\subsubsection{Proper push-forward}(\cite[2.2.7 (2)]{DJK})
\label{num:DJK2272}
If $p:X\to Y$ is a proper morphism and $v$ is a virtual vector bundle on $Y$, then there is a push-forward map $p_*:H(X/S,p^*v)\to H(Y/S,v)$.

\subsubsection{(Refined) pullback}(\cite[Def. 4.2.5]{DJK})
\label{num:DJK425}
Consider a Cartesian square as in~\eqref{eq:Cartsqdelta}, and assume that $p$ is a local complete intersection (\ref{num:lci}). Let $v$ be a virtual vector bundle on $X$. Then there is a map
\begin{align}
\label{eq:refgysin}
\Delta^!:H(X/S,v)\to H(Y/S,g^*\tau_p+v_{|Y})
\end{align}
induced by the \textbf{refined fundamental class} $\Delta^*\eta_p$. In particular, if $q:Y\to X$ is \'etale, then we have $q^*:H(X/S,v)\to H(Y/S,v_{|Y})$ (\cite[2.2.7 (3)]{DJK}).

\begin{lemma}
\label{lm:delta*!}
Consider a Cartesian square as in~\eqref{eq:Cartsqdelta}
Assume that $p$ is a regular closed immersion and $\mathbbold{1}_S$ is $p$-pure. Let $v$ be a virtual vector bundle on $X$.
Then the following diagram is commutative:
\begin{align}
\begin{split}
  \xymatrix{
    H(X/S,v) \ar[r]^-{\Delta'^*} \ar[d]_-{\Delta^!} & H(Y/Z,v_{|Y}) \\
    H(Y/S,v_{|Y}-g^*N_ZS) \ar[ru]_-{\sim} & 
  }
\end{split}
\end{align}
\proof
This follows from the definition of the map $\Delta^!$ in~\eqref{eq:refgysin} and \cite[2.2.13]{DJK}.
\endproof

\end{lemma}

\subsubsection{Specialization}(\cite[4.5.6]{DJK})
\label{num:spbiv}
Let $i:Z\to S$ be a regular closed immersion and let $j:U\to S$ be the open complement. Denote by $e(N_ZS)$ the Euler class of the normal bundle of $Z$ in $S$, and assume that there is a null-homotopy $e(N_ZS)\simeq0$. Let $f:X\to S$ be a separated morphism of finite type, and consider the Cartesian square
\begin{align}
\begin{split}
  \xymatrix@=11pt{
    X_Z \ar[r]^-{i_X} \ar[d]_-{f_Z} & X \ar[d]^-{f} & X_U \ar[d]^-{f_U} \ar[l]_-{j_X}\\
    Z \ar[r]^-{i} & S & U  \ar[l]_-{j}.
  }
\end{split}
\end{align}
For any object $A\in\mathcal{SH}(X)$, we have the composition
\begin{align}
\label{eq:DJKspint}
i_{X*}i_X^!A\to A\to i_{X*}i_X^*A\to i_{X*}(i_X^!A\otimes Th(-N_ZS)_{|X_Z})
\end{align}
where the last map is induced by the refined fundamental class of $i$ (\cite[Def. 4.2.5]{DJK}). By the self-intersection formula (\cite[Cor. 4.2.3]{DJK}), the map~\eqref{eq:DJKspint} agrees with the multiplication by the class $f_Z^*e(N_ZS)$, which is null-homotopic by assumption. By the localization sequence, we obtain a natural transformation of functors
\begin{align}
\label{eq:DJK456a}
j_{X*}j_X^!A\to i_{X*}(i_X^!A\otimes Th(-N_ZS)_{|X_Z}).
\end{align}
Now assume that $\mathbbold{1}_S$ is $i$-pure in the sense of~\ref{num:lci}, and let $v$ be a virtual vector bundle on $X$. Then the map~\eqref{eq:DJK456a} induces the following \textbf{specialization map}:
\begin{align}
\label{eq:DJK456b}
H(X_U/U,v_{|X_U})=H(X_U/S,v_{|X_U})\to H(X_Z/S,v_{|X_Z}-(N_ZS)_{|X_Z})\simeq H(X_Z/Z,v_{|X_Z}).
\end{align}
\begin{lemma}
\label{lm:toptriv}
Assume that $i:Z\to S$ has a smooth retraction $S\to Z$.
Let $Y$ be a $Z$-scheme and let $v$ be a virtual vector bundle on $Y$. Let $X\to Y\times_ZS$ be a morphism inducing an isomorphism on the reduced schemes $X_{red}\to (Y\times_ZS)_{red}$.
Then the following composition is an isomorphism:
\begin{align}
H(Y/Z,v)
\xrightarrow{\eqref{eq:horBC}}
H(Y\times_ZU/U,v_{|Y\times_ZU})
\simeq 
H(X_U/U,v_{|X_U})
\xrightarrow{\eqref{eq:DJK456b}}
H(X_Z/Z,v_{|X_Z}).
\end{align}

\end{lemma}
\proof
It follows from the localization sequence that for any $T$-scheme $W$ and any virtual vector bundle $w$ on $W$, the canonical morphism $W_{red}\to W$ induces an isomorphism $H(W_{red}/T,w_{|W_{red}})\simeq H(W/T,w)$. Therefore to prove the claim we may assume that $X=Y\times_ZS$ is a trivial family. By construction of the map~\eqref{eq:DJK456b}, it suffices to show that the composition
\begin{align}
\label{eq:delta12*!}
H(Y/Z,v)
\xrightarrow{\Delta_2'^*}
H(Y\times_ZS/S,v_{|Y\times_ZS})
\xrightarrow{\Delta_1^!}
H(Y/S,v-(N_ZS)_{|Y})
\simeq
H(Y/Z,v)
\end{align}
is identity, where $\Delta_1$ and $\Delta_2$ are Cartesian squares in the diagram
\begin{align}
\begin{split}
  \xymatrix@=11pt{
    Y \ar[r]^-{} \ar[d]_-{} \ar@{}[rd]|{\Delta_1} & Y\times_ZS \ar[d]^-{} \ar[r]_-{} \ar@{}[rd]|{\Delta_2} & Y \ar[d]^-{} \\
    Z \ar[r]^-{i} & S \ar[r]_-{} & Z.
  }
\end{split}
\end{align}
By Lemma~\ref{lm:delta*!}, the map~\eqref{eq:delta12*!} agrees with the composition
\begin{align}
H(Y/Z,v)
\xrightarrow{\Delta_2'^*}
H(Y\times_ZS/S,v_{|Y\times_ZS})
\xrightarrow{\Delta_1'^*}
H(Y/Z,v)
\end{align}
which is identity since the formation of the map~\eqref{eq:horBC} is compatible with composition of squares.
\endproof

\subsection{On the generators of motivic categories}
\subsubsection{}
Recall the following de Jong-Gabber alteration result:
\begin{lemma}[\textrm{\cite[X Thm. 2.1, 3.5]{ILO}}]
\label{lm:gabber}
Let $k$ be a field, $l$ a prime number different from the characteristic of $k$, $X$ an integral separated $k$-scheme finite type, $Z\subset X$ a nowhere dense closed subset. Then there exist
\begin{itemize}
\item
a finite purely inseparable field extension $k'$ of $k$ of degree prime to $l$
\item
a projective, surjective, generically finite, maximally dominating morphism $h:X'\to X$
\end{itemize}
such that
\begin{enumerate}
\item
$X'$ is integral, smooth and quasi-projective over $k'$;
\item
the generic degree of $h$ is prime to $l$;
\item
\label{num:gabbergalois}
the morphism between the function fields of $X'$ and $X$ induced by $h$ is a finite Galois field extension followed by a finite purely inseparable field extension;
\item
$h^{-1}(Z)$ is the support of a strict normal crossings divisor.
\end{enumerate}
\end{lemma}
\proof
While the other points have been proved in \cite[X Thm. 2.1]{ILO}, let us explain how to achieve the point~\eqref{num:gabbergalois} using the arguments in the proof of \cite[X Thm. 3.5]{ILO}. One uses the terminology ``separable $l'$-alteration'' as in loc. cit.. Denote by $k^{perf}$ the perfection of the field $k$. Since $X\times_kk^{perf}$ is separated of finite type over the perfect field $k^{perf}$, by \cite[X Thm. 3.5]{ILO} one can find a separable $l'$-alteration $X_1\to X\times_kk^{perf}$ with $X_1$ smooth over $k^{perf}$. Since $k^{perf}$ is the colimit of finite purely inseparable extensions of $k$, by \cite[Thm. 8.10.5]{EGA43} and \cite[Prop. 17.7.8]{EGA4} there exists a finite purely inseparable extension $k'$ of $k$ and a separable $l'$-alteration $X'\to X\times_{k}k'$ such that $X_1\simeq X'\times_{k'}k^{perf}$, and $X'$ is smooth over $k'$. Similarly we may obtain that the inverse image of $Z$ is the support of a strict normal corssings divisor. Then the composition $h:X'\to X\times_{k}k'\to X$ gives us what we need.
\endproof

\subsubsection{}
Recall that a morphism of schemes is a \emph{universal homeomorphism} if and only if it is integral, radicial and surjective (\cite[Cor. 18.12.11]{EGA4}). We have the following topological invariance of $\mathcal{SH}$ due to Elmanto-Khan:
\begin{lemma}[\textrm{\cite[Thm. 2.1.1]{EK}}]
\label{lm:EK211}
For $\mathfrak{P}$ a collection of prime numbers, the motivic $\infty$-category $\mathcal{SH}(-)[\mathfrak{P}^{-1}]$ is \emph{topologically invariant} when restricted to schemes where every prime outside $\mathfrak{S}$ is invertible, that is, for any universal homeomorphism of such schemes $f:T\to S$, the functor $f^*:\mathcal{SH}(S)[\mathfrak{P}^{-1}]\to\mathcal{SH}(T)[\mathfrak{P}^{-1}]$ is an equivalence.

In particular, for any prime $p$, the motivic $\infty$-category 
$\mathcal{SH}(-)[1/p]$ is topologically invariant when restricted to schemes of characteristic $p$.
\end{lemma}

\subsubsection{}
For any morphism $g:T\to S$, we denote $M^{BM}(T/S)=g_!\mathbbold{1}_T\in\mathcal{SH}(S)$. If $W$ is a closed subscheme of $T$, then by localization there is a canonical cofiber sequence
\begin{align}
\label{eq:BMloc}
M^{BM}(T-W/S)\to M^{BM}(T/S)\to M^{BM}(W/S).
\end{align}
The following lemma slightly strengthens \cite[Prop. 3.1.3]{EK}:
\begin{lemma}
\label{lm:CD7.2}
Let $k$ be a field of exponential characteristic $p$, and let $\mathcal{T}$ be a $\mathbb{S}[1/p]$-linear motivic $\infty$-category. 
Let $S$ be a $k$-scheme of finite type and let $Z$ be a closed subset of $S$. 
Then the subcategory of constructible objects $\mathcal{T}_c(S)$ is generated as a thick subcategory by elements of the form $f_*\mathbbold{1}_X(n)$, where $f:X\to S$ is a projective morphism such that $X$ is smooth over a finite purely inseparable extension of $k$, and $f^{-1}(Z)$ is either empty, or the whole $X$, or the support of a strict normal crossing divisor, and $n$ is an integer.
\end{lemma}
\proof
By universality of $\mathcal{SH}$ (\cite[Thm. 7.14]{DG}), we may assume that $\mathcal{T}=\mathcal{SH}(-)[1/p]$. 
We follow the lines of \cite[Prop. 7.2]{CD1}. Let $\mathcal{G}$ be the family of objects of the form $f_*\mathbbold{1}_X(n)$, where $f:X\to S$ is a projective morphism such that $X$ is smooth over a finite purely inseparable extension of $k$ and $f^{-1}(Z)$ is either empty, or the whole $X$, or the support of a strict normal crossing divisor, and $n$ is an integer. Let $l$ be a prime number different from $p$, and let $\mathcal{C}$ be the thick subcategory of $\mathcal{SH}(S,\mathbb{Z}_{(l)})$ generated by the elements of $\mathcal{G}$. We know that the elements of $\mathcal{C}$ are constructible, and it suffices to show that $\mathcal{C}=\mathcal{SH}_c(S,\mathbb{Z}_{(l)})$. 

By \cite[Lemme 2.2.23]{Ayo}, $\mathcal{SH}_c(S,\mathbb{Z}_{(l)})$ is the thick subcategory of $\mathcal{SH}(S,\mathbb{Z}_{(l)})$ generated by the elements of the form $f_*\mathbbold{1}_X(n)$, where $f:X\to S$ is a projective morphism and $n$ is an integer. We claim that $f_*\mathbbold{1}_X(n)$ lies in $\mathcal{C}$ for all such morphisms $f$ and integer $n$, and we use induction on the dimension of $X$. By $cdh$ descent (\cite[Prop. 3.7]{Cis1}), we may reduce to the case where $X$ is integral. If $X$ is empty, there is nothing to prove. So by induction we may assume that the claim holds for all projective $S$-schemes $X_1$ whose dimension is smaller than $X$. As in Lemma~\ref{lm:gabber}, let $h:X'\to X$ be a projective, surjective, generically finite morphism of degree prime to $l$, such that $X'$ is integral and smooth over a finite purely inseparable extension of $k$ and $h^{-1}(Z)$ is either empty, or the whole $X$, or the support of a strict normal crossing divisor, and such that the finite field extension induced by $h$ on generic points $\widetilde{h}:\operatorname{Spec}(L)\to\operatorname{Spec}(K)$ is a finite Galois extension followed by a finite purely inseparable extension. 
Then by \cite[Lemmas B.3 and B.4]{LYZR} and Lemma~\ref{lm:EK211}, the canonical map $\widetilde{h}_!\mathbbold{1}_L\to\mathbbold{1}_K$ has a section in $\mathcal{SH}_c(K,\mathbb{Z}_{(l)})$. By continuity of $\mathcal{SH}$ as explained in the proof of \cite[Lemma 2.4.6]{BD}, there exists a non-empty open subscheme $U$ of $X$ such that for $V=h^{-1}(U)$, $M^{BM}(U/S)$ is a direct summand of $M^{BM}(V/S)$ in $\mathcal{SH}(S,\mathbb{Z}_{(l)})$. On the other hand, by the localization sequence~\eqref{eq:BMloc}, for any non-empty open subscheme $V$ of $X'$, the object $M^{BM}(V/S)$ lies in $\mathcal{C}$, since $X'-V$ has dimension smaller than $X$. It follows that $M^{BM}(U/S)$ also lies in $\mathcal{C}$. We conclude by using the localization sequence~\eqref{eq:BMloc} and the induction assumption again.
\endproof

\begin{remark}
Lemma~\ref{lm:CD7.2} has several variants and predecessors in the literature: see \cite[Prop. 2.2.27]{Ayo}, \cite[Prop. 7.2]{CD1}, \cite[Cor. 2.4.8]{BD}, \cite[Prop. 3.1.3]{EK}.
\end{remark}

\section{Correspondences and the trace}
\label{sec:trace}
\subsection{Relative K\"unneth formulas}
\label{sec:relkun}
In this section we prove K\"unneth formulas over a base scheme. We use the techniques developed in \cite[\S2]{JY} and follow the ideas in \cite[\S2.1]{LZ}.
\begin{definition}[\textrm{\cite[Def. 2.1.7]{JY}}]
\label{def:locacy}
Let $f:X\to S$ be a morphism of schemes and $K\in \mathcal{SH}(X)$. We say that $K$ is \textbf{strongly locally acyclic} over $S$ if for any Cartesian square
\begin{align}
\begin{split}
  \xymatrix@=10pt{
    Y \ar[r]^-{q} \ar[d]_-{g} & X \ar[d]^-{f}\\
    T \ar[r]_-{p} & S
  }
\end{split}
\end{align}
and any object $L\in\mathcal{SH}(T)$, the canonical map $K\otimes f^*p_*L\to q_*(q^*K\otimes g^*L)$ is an isomorphism. We say that $K$ is \textbf{universally strongly locally acyclic} (abbreviated as \textbf{USLA}) over $S$ if for any morphism $T\to S$, the base change $K_{|X\times_ST}$ is strongly locally acyclic over $T$.

\end{definition}

The USLA objects form a localizing subcategory. Typical examples are given by dualizable objects for smooth morphisms:
\begin{lemma}
\label{lm:smUSLA}
If $f:X\to S$ is smooth, any dualizable object in $\mathcal{SH}(X)$ is USLA over $S$.
\end{lemma}
\proof
This follows from smooth base change and \cite[Prop. 3.2]{FHM}.
\endproof
Another important case is the one where $S$ is the spectrum of a field, in which, assuming resolution of singularities, every object is USLA:
\begin{lemma}[\textrm{\cite[Cor. 2.1.14]{JY}}]
\label{lm:USLA}
Let $k$ be a field and let $\mathfrak{P}$ be a collection of primes such that at least one of the following conditions is satisfied:
\begin{enumerate}
\item $k$ is a perfect field which satisfies strong resolution of singularities (\ref{num:RS});

\item 
\label{num:invertp}
$\mathfrak{P}$ contains $p$, the exponential characteristic of $k$.
\end{enumerate}
Then for every separated $k$-scheme of finite type $X$, every object of $\mathcal{SH}(X)[\mathfrak{P}^{-1}]$ is USLA over $k$.
\end{lemma}

\begin{lemma}
\label{lm:LZ2.22.3}
Let $Y\to S$ be a morphism and let $M\in\mathcal{SH}(Y)$ be USLA over $S$. Let $f:X\to X'$ be a morphism of $S$-schemes, and let $f_Y=f\times_Sid_Y:X\times_SY\to X'\times_SY$.
\begin{enumerate}
\item \label{num:LZ2.2}
For any $L\in\mathcal{SH}(X)$, there is a canonical isomorphism
\begin{align}
f_*L\boxtimes_SM\simeq f_{Y*}(L\boxtimes_SM).
\end{align}
\item
For any $L'\in\mathcal{SH}(X')$, there is a canonical isomorphism
\begin{align}
\label{eq:LZ2.3}
f^!L'\boxtimes_SM\simeq f_Y^!(L'\boxtimes_SM).
\end{align}
\end{enumerate}
\end{lemma}
\proof
\begin{enumerate}
\item This is a reformulation of the definition of being USLA.
\item We may assume that $f$ is smooth or a closed immersion. The smooth case follows from relative purity. If $f$ is a closed immersion, let $j$ be the complementary open immersion and let $j_Y=j\times_Sid_Y$. Then we have a commutative diagram
\begin{align}
\begin{gathered}
  \xymatrix{
    f^!L'\boxtimes_SM \ar[r]^-{} \ar[d]_-{} & f^*L'\boxtimes_SM \ar[d]^-{} \ar[r]^-{} & f^*j_*j^*L'\boxtimes_SM \ar[r]^-{+1} \ar[d]^-{} &\\
    f_Y^!(L'\boxtimes_SM) \ar[r]^-{} & f_Y^*(L'\boxtimes_SM) \ar[r]^-{} & f_Y^*j_{Y*}j_Y^*(L'\boxtimes_SM) \ar[r]^-{+1} &
  }
\end{gathered}
\end{align}
where both rows are distinguished triangles. The middle vertical map is an isomorphism, and the right vertical map is an isomorphism by~\eqref{num:LZ2.2}. Therefore by five lemma the left vertical map is also an isomorphism.
\end{enumerate}
\endproof
The following corollary is the special case of~\eqref{eq:LZ2.3} for $L'=\mathbbold{1}$:
\begin{corollary}
\label{cor:LZ2.4}
Let $X$ and $Y$ be two $S$-schemes and denote by $p_Y:X\times_SY\to Y$ the projection. Let $M\in\mathcal{SH}(Y)$ be USLA over $S$. Then there is a canonical isomorphism
\begin{align}
\mathcal{K}_{X/S}\boxtimes_SM\simeq p_Y^!M
\end{align}
where $p_Y:X\times_SY\to Y$ is the projection.
\end{corollary}

\begin{proposition}
\label{prop:LZ2.5}
Let $X$ and $Y$ be two $S$-schemes and denote by $p_X:X\times_SY\to X$ and $p_Y:X\times_SY\to Y$ the projections. Let $M\in\mathcal{SH}(Y)$ be USLA over $S$ and let $L\in\mathcal{SH}_c(X)$ be a constructible object. Then there is a canonical isomorphism
\begin{align}
\label{eq:Kunrel}
\mathbb{D}_{X/S}(L)\boxtimes_SM\simeq\underline{Hom}(p_X^*L,p_Y^!M).
\end{align}
\end{proposition}

\proof
We may assume that $L=\phi_*\mathbbold{1}_W$ for some proper morphism $\phi:W\to X$. Denote by $\phi_Y:W\times_SY\to X\times_SY$ the base change. Then we have 
\begin{align}
\begin{split}
\phi_*\mathcal{K}_{W/S}\boxtimes_SM
&\simeq
\phi_{Y*}(\mathcal{K}_{W/S}\boxtimes_SM)
\simeq
\phi_{Y*}\underline{Hom}(\mathbbold{1}_{W\times_SY},\phi_{Y}^!p_Y^!M)\\
&\simeq
\underline{Hom}(p_X^*\phi_*\mathbbold{1}_W,p_Y^!M)
\end{split}
\end{align}
where the first isomorphism follows from Lemma~\ref{lm:LZ2.22.3} and the second isomorphism from Corollary~\ref{cor:LZ2.4}.
\endproof

\begin{corollary}
\label{cor:LZ211}
Let $X$ be an $S$-scheme and let $L\in\mathcal{SH}_c(X)$ be a constructible object which is USLA over $S$. Then the canonical map $L\to \mathbb{D}_{X/S}(\mathbb{D}_{X/S}(L))$ is an isomorphism.
\end{corollary}
\proof
This can be proved using the fact that $L$ is a \emph{dualizable object} in the \emph{category of cohomological correspondences}, which is a consequence of Proposition~\ref{prop:LZ2.5}. We refer to \cite[Prop. 2.11]{LZ} for more details.
\endproof

\subsection{Correspondences and the trace}
\label{sec:corrtr}
Using Proposition~\ref{prop:LZ2.5}, we generalize the trace map defined in \cite[Prop. 3.2.8]{JY} to more general base schemes.

\subsubsection{}
For $X\to S$ a morphism, denote by $p_1,p_2:X\times_SX\to X$ the projections. A \textbf{correspondence} is a morphism of the form $c:C\to X\times_SX$. We denote by $c_1,c_2:C\to X$ the compositions of $c$ with $p_1$ and $p_2$. Note that every $S$-endormorphism $f:X\to X$ can be viewed as a correspondence via the transpose of the graph morphism $X\xrightarrow{(f, id)}X\times_SX$.

\subsubsection{}
Consider the following Cartesian diagram
\begin{align}
\begin{split}
  \xymatrix@=10pt{
    Fix(c) \ar[r]^-{c'} \ar[d]_-{\delta'} & X \ar[d]^-{\delta_{X/S}}\\
    C \ar[r]^-{c} & X\times_SX.
  }
\end{split}
\end{align}
Let $K\in\mathcal{SH}_c(X)$ be USLA over $S$, and let $v$ be a virtual vector bundle on $C$. A \textbf{(cohomological) correspondence over $c$} is a map of the form $u:c_1^*K\to c_2^!K\otimes Th(v)$. Given such a map, we have the composition
\begin{align}
\label{eq:JY3281}
\begin{split}
Th(-v)
\xrightarrow{u}
\underline{Hom}(c_1^*K,c_2^!K)
\overset{}{\simeq}
c^!\underline{Hom}(p_1^*K,p_2^!K)
\overset{\eqref{eq:Kunrel}}{\simeq}
c^!(\mathbb{D}_{X/S}(K)\boxtimes_S K).
\end{split}
\end{align}
which gives rise to the following map
\begin{align}
\label{eq:gen_tr}
\begin{split}
c'_!Th(-v_{|Fix(c)})
&\simeq
c'_!\delta'^*Th(-v)
\simeq
\delta_{X/S}^*c_!Th(-v)
\xrightarrow{\eqref{eq:JY3281}}
\delta_{X/S}^*c_!c^!(\mathbb{D}_{X/S}(K)\boxtimes_SK)\\
&\to
\delta_{X/S}^*(\mathbb{D}_{X/S}(K)\boxtimes_SK)
=
\mathbb{D}_{X/S}(K)\otimes K
\simeq
K\otimes\mathbb{D}_{X/S}(K)
\to
\mathcal{K}_{X/S}.
\end{split}
\end{align}

\begin{definition}
\label{def:tracerel}
Let $K\in\mathcal{SH}_c(X)$ be USLA over $S$.
For a correspondence $u:c_1^*K\to c_2^!K\otimes Th(v)$, we denote by 
\begin{align}
\label{eq:traceu}
Tr(u/S):Th(-v_{|Fix(c)})\to \mathcal{K}_{Fix(c)/S}
\end{align}
the map obtained by adjunction from the map~\eqref{eq:gen_tr}, called the \textbf{trace map}. This construction gives rise to a canonical map 
\begin{align}
\label{eq:ustrace}
Map(c_1^*K,c_2^!K\otimes Th(v))\xrightarrow{Tr(-/S)}H(Fix(c)/S,-v_{|Fix(c)}).
\end{align}
\end{definition}

\subsubsection{}
The following lemma states the additivity of traces along distinguished triangles, where the structure of higher categories plays a key role. We refer to \cite[\S4]{JY} for a detailed discussion.
\begin{lemma}[Additivity of traces]
\label{lm:additivity}
Let $L\to M\to N$ be a cofiber sequence in $\mathcal{SH}_c(X)$ of USLA objects over $S$, and let
\begin{align}
\begin{split}
  \xymatrix@=10pt{
    c_1^*L \ar[r]^-{} \ar[d]_-{u_L} & c_1^*M \ar[r]^-{} \ar[d]^-{u_M} & c_1^*N \ar[d]^-{u_N}\\ 
    c_2^!L\otimes Th(v) \ar[r]^-{} & c_2^!M\otimes Th(v) \ar[r]^-{} & c_2^!N\otimes Th(v)
  }
\end{split}
\end{align}
be a morphism of cofiber sequences (in the $\infty$-categorical sense).
Then there is a canonical homotopy between $Tr(u_M/S)$ and $Tr(u_L/S)+Tr(u_N/S)$ as maps $Th(-v_{|Fix(c)})\to \mathcal{K}_{Fix(c)/S}$.
\end{lemma}
\proof
The additivity is proved in \cite[Prop. 4.2.6]{JY} when $S$ is the spectrum of a field, using the language of motivic derivators (\cite[Def. 2.4.48]{Ayo}); by a similar argument the proof also works for a general base scheme $S$: the only technical difficulty lies in \cite[Lemma 4.1.15]{JY} where, even though the local duality fails to hold in general, the proof can be completed by applying Corollary~\ref{cor:LZ211}. 
\endproof

\begin{remark}
\label{rm:verdiertrace}
We may readily extend the trace map to the Verdier pairing over a base $\langle u,v\rangle_S$, in the same fashion as \cite[III 4.1]{SGA5}, \cite[Def. 3.1.8]{JY} or \cite[\S2.4]{LZ}. It turns out that the computation of the Verdier pairing reduces to that of the trace map via the identity $\langle u,v\rangle_S=\langle vu,1\rangle_S$, see \cite[Prop. 3.2.5]{JY}, and the additivity in Lemma~\ref{lm:additivity} can also be extended as in \cite[Thm. 4.2.8]{JY}. In this paper we will not discuss these variants, and we will focus on trace maps.
\end{remark}

\section{Operations on correspondences}
\label{sec:opcorr}
We discuss three types of operations on correspondences: base change, pullback and push-forward, which are all compatible with the trace map. 

\subsection{Base change}
\label{sec:bc}
\subsubsection{}
Let $c:C\to X\times_SX$ and $T\to S$ be two morphisms. Let $Y=X\times_ST$. Then there is a canonical morphism $c_T:C_T=C\times_ST\to Y\times_TY$. Let $K\in\mathcal{SH}(X)$ and let $v$ be a virtual vector bundle on $C$. Given a correspondence $u:c_1^*K\to c_2^!K\otimes Th(v)$, we have the following composition
\begin{align}
c_{T1}^*K_{|Y}=(c_1^*K)_{|C_T}\xrightarrow{u}(c_2^!K\otimes Th(v))_{|C_T}\to c_{T2}^!K_{|Y}\otimes Th(v_{|C_T}).
\end{align}
This construction gives rise to a canonical map
\begin{align}
\label{eq:BCcorr}
Map(c_1^*K, c_2^!K\otimes Th(v))\to Map(c_{T1}^*K_{|Y}, c_{T2}^!K_{|Y}\otimes Th(v_{|C_T}))
\end{align}
called the \textbf{base change of correspondences} from $S$ to $T$. On the other hand, we have a canonical Cartesian square
\begin{align}
\begin{split}
  \xymatrix@=10pt{
    Fix(c_T) \ar[r]^-{} \ar[d]_-{} \ar@{}[rd]|{\Delta} & T \ar[d]^-{}\\
    Fix(c) \ar[r]^-{} & S
  }
\end{split}
\end{align}
so by~\ref{num:DJK2271} there is a canonical map
\begin{align}
\label{eq:BCbiv}
\Delta^*:H(Fix(c)/S)\to H(Fix(c_T)/T).
\end{align}
The following lemma states that base change maps are compatible with the trace map, which generalizes \cite[Prop. 6.2.16]{JY}:
\begin{lemma}
\label{lm:BCtrace}
Let $K\in\mathcal{SH}_c(X)$ be USLA over $S$. Then the following diagram is commutative:
\begin{align}
\begin{split}
  \xymatrix{
    Map(c_1^*K, c_2^!K\otimes Th(v)) \ar[r]^-{Tr(-/S)} \ar[d]_-{\eqref{eq:BCcorr}} & H(Fix(c)/S,-v_{|Fix(c)}) \ar[d]^-{\eqref{eq:BCbiv}}\\
    Map(c_{T1}^*K_{|Y}, c_{T2}^!K_{|Y}\otimes Th(v_{|C_T})) \ar[r]^-{Tr(-/T)} & H(Fix(c_T)/T,-v_{|Fix(c_T)}).
  }
\end{split}
\end{align}

\end{lemma}
The proof is very similar to that of \cite[Prop. 6.2.16]{JY}, which leave as exercise.

\subsection{Push-forward} 

\subsubsection{}
\label{num:corrpf}
Let $f:X\to Y$ be a morphism of $S$-schemes. Let $c:C\to X\times_SX$ and $d:D\to Y\times_SY$ be two correspondences. Let $p:C\to D$ be a morphism such that the following diagram is commutative:
\begin{align}
\label{eq:Var1.1big}
\begin{split}
  \xymatrix@=10pt{
    C \ar[r]^-{c} \ar[d]_-{p} \ar@{}[rd]|{\Delta} & X\times_SX \ar[d]^-{f\times_Sf}\\
    D \ar[r]_-{d} & Y\times_SY.
  }
\end{split}
\end{align}
Assume that one of the following conditions hold:
\begin{enumerate}
\item The following commutative square is Cartesian:
\begin{align}
\label{eq:Var1.1left}
\begin{split}
  \xymatrix@=10pt{
    C \ar[r]^-{c_1} \ar[d]_-{p} & X \ar[d]^-{f}\\
    D \ar[r]^-{d_1} & Y;
  }
\end{split}
\end{align}
\item Both morphisms $p$ and $f$ are proper.
\end{enumerate}
Then there is a natural transformation of functors $d_1^*f_!\to p_!c_1^*$. Let $v$ be a virtual vector bundle on $D$. Then for any $K\in\mathcal{SH}(X)$ and any correspondence $u:c_1^*K\to c_2^!K\otimes Th(v_{|C})$, we have the following composition
\begin{align}
d_1^*f_!K\to p_!c_1^*K\xrightarrow{u}p_!(c_2^!K\otimes Th(v_{|C}))\simeq p_!c_2^!K\otimes Th(v)\to d_2^!f_!K\otimes Th(v).
\end{align}
This construction gives rise to a canonical map
\begin{align}
\label{eq:pfcorr}
\Delta_!:
Map(c_1^*K, c_2^!K\otimes Th(v_{|C}))\to Map(d_1^*f_!K, d_2^!f_!K\otimes Th(v))
\end{align}
called the \textbf{push-forward of correspondences}.

\subsubsection{}
Assume that both $p$ and $f$ are proper. We have a canonical commutative square
\begin{align}
\begin{split}
  \xymatrix@=10pt{
    Fix(c) \ar[r]^-{c'} \ar[d]_-{q} & X \ar[d]^-{}\\
    Fix(d) \ar[r]^-{d'} & Y
  }
\end{split}
\end{align}
with $q$ proper, which by~\ref{num:DJK2272} induces a canonical map
\begin{align}
\label{eq:pfbiv}
q_*:H(Fix(c)/S,-v_{|Fix(c)})\to H(Fix(d)/S,-v_{|Fix(d)}).
\end{align}
\begin{lemma}
\label{lm:pfUSLA}
Let $f:X\to Y$ be a proper morphism of $S$-schemes. For any object $K\in\mathcal{SH}(X)$ USLA over $S$, the object $f_*K$ is also USLA over $S$.
\end{lemma}
\proof
For any Cartesian square
\begin{align}
\begin{split}
  \xymatrix@=10pt{
    X_T \ar[d]_-{t} \ar[r]^-{g} &  Y_T \ar[d]^-{s} \ar[r]^-{q} & T \ar[d]^-{r}\\
    X  \ar[r]^-{f}  & Y \ar[r]^-{p} & S
  }
\end{split}
\end{align}
and any object $L\in\mathcal{SH}(T)$, we have canonical isomorphisms
\begin{align}
\begin{split}
f_*K\otimes p^*r_*L
&\simeq 
f_*(K\otimes f^*p^*r_*L)
\simeq
f_*t_*(t^*K\otimes g^*q^*L)
=
s_*g_*(t^*K\otimes g^*q^*L)\\
&\simeq
s_*(g_*t^*K\otimes q^*L)
\simeq
s_*(s^*f_*K\otimes q^*L)
\end{split}
\end{align}
where we use the properness of $f$ and $g$ and the fact that $K$ is USLA over $S$. The same property also holds after any base change, which implies that $f_*K$ is USLA over $S$.
\endproof
The following lemma states that proper push-forwards are compatible with the trace map:
\begin{lemma}
\label{lm:LZ221}
Let $K\in\mathcal{SH}_c(X)$ be USLA over $S$. Then the following diagram is commutative:
\begin{align}
\begin{split}
  \xymatrix{
    Map(c_1^*K, c_2^!K\otimes Th(v_{|C})) \ar[r]^-{Tr(-/S)} \ar[d]_-{\eqref{eq:pfcorr}} & H(Fix(c)/S,-v_{|Fix(c)}) \ar[d]^-{\eqref{eq:pfbiv}}\\
    Map(d_1^*f_*K, d_2^!f_*K\otimes Th(v)) \ar[r]^-{Tr(-/S)} & H(Fix(d)/S,-v_{|Fix(d)}).
  }
\end{split}
\end{align}
where the lower horizontal map is well-defined by Lemma~\ref{lm:pfUSLA}.
\end{lemma}
Lemma~\ref{lm:LZ221} is a particular case of the Lefschetz-Verdier formula. For details of the proof, see \cite[III 4.4]{SGA5}, \cite[Prop. 1.2.5]{Var}, \cite[Thm. 3.3.2]{YZ}, \cite[Thm. 3.2.18]{Cis} or \cite[Thm. 2.20]{LZ}.

\begin{remark}
Let $p:X\to S$ be a smooth proper morphism. By Lemma~\ref{lm:smUSLA}, the object $\mathbbold{1}_X$ is USLA over $S$. In diagram~\eqref{eq:Var1.1big}, consider the case $C=X$, $D=Y=S$, $f=p$, $c=\delta_{X/S}$, $d=id_S$ and $v=0$. Then Lemma~\ref{lm:LZ221} applied to the identity map recovers the motivic Gauss-Bonnet formula (\cite[Thm. 4.6.1]{DJK}, \cite[Thm. 5.3]{LR}), see also \cite[Rmk. 5.1.11 (2)]{JY}.
\end{remark}

\subsection{Pullback}

\subsubsection{}
\label{num:corrpb}
Consider the situation of Diagram~\eqref{eq:Var1.1big}.
Assume that one of the following conditions hold:
\begin{enumerate}
\item The commutative square 
\begin{align}
\label{eq:Var1.1right}
\begin{split}
  \xymatrix@=11pt{
    C \ar[r]^-{c_2} \ar[d]_-{p} \ar@{}[rd]|{\Delta} & X \ar[d]^-{f}\\
    D \ar[r]^-{d_2} & Y,
  }
\end{split}
\end{align}
induces an isomorphism $C_{red}\to (D\times_YX)_{red}$; 
\item Both morphisms $p$ and $f$ are \'etale.
\end{enumerate}
Then there is a natural transformation of functors $p^*d_2^!\to c_2^!f^*$. Let $v$ be a virtual vector bundle on $D$. Then for any $K\in\mathcal{SH}(Y)$ and any correspondence $u:d_1^*K\to d_2^!K\otimes Th(v)$, we have the following composition
\begin{align}
c_1^*f^*K=p^*d_1^*K\xrightarrow{u}p^*(d_2^!K\otimes Th(v))\to c_2^!f^*K\otimes Th(v_{|C}).
\end{align}
This construction gives rise to a canonical map
\begin{align}
\label{eq:pbcorr}
(-)_{|\Delta}:Map(d_1^*K, d_2^!K\otimes Th(v))\to Map(c_1^*f^*K, c_2^!f^*K\otimes Th(v_{|C}))
\end{align}
called the \textbf{pull-back of correspondences}.

\subsubsection{}
Assume that both $p$ and $f$ are \'etale. We have a canonical commutative square
\begin{align}
\begin{split}
  \xymatrix@=10pt{
    Fix(c) \ar[r]^-{c'} \ar[d]_-{q} & X \ar[d]^-{}\\
    Fix(d) \ar[r]^-{d'} & Y
  }
\end{split}
\end{align}
with $q$ \'etale, which by~\ref{num:DJK425} induces a canonical map
\begin{align}
\label{eq:pbbiv}
q^*:H(Fix(d)/S,-v_{|Fix(d)})\to H(Fix(c)/S,-v_{|Fix(c)}).
\end{align}

\subsubsection{}
The proof of the following Lemma is very similar to Lemma~\ref{lm:pfUSLA} and is left as an exercise:
\begin{lemma}
\label{lm:pbUSLA}
Let $f:X\to Y$ be a smooth morphism of $S$-schemes. For any object $K\in\mathcal{SH}(Y)$ USLA over $S$, the object $f^*K$ is also USLA over $S$.
\end{lemma} 

\subsubsection{}
The following lemma states that \'etale pullbacks are compatible with the trace map:
\begin{lemma}
\label{lm:traceetpb}
Let $K\in\mathcal{SH}_c(Y)$ be USLA over $S$. Then the following diagram is commutative:
\begin{align}
\begin{split}
  \xymatrix{
    Map(d_1^*K, d_2^!K\otimes Th(v)) \ar[r]^-{Tr(-/S)} \ar[d]_-{\eqref{eq:pbcorr}} & H(Fix(d)/S,-v_{|Fix(d)}) \ar[d]^-{\eqref{eq:pbbiv}}\\
    Map(c_1^*f^*K, c_2^!f^*K\otimes Th(v_{|C})) \ar[r]^-{Tr(-/S)} & H(Fix(c)/S,-v_{|Fix(c)}).
  }
\end{split}
\end{align}
where the lower horizontal map is well-defined by Lemma~\ref{lm:pbUSLA}.
\end{lemma}
The proof of Lemma~\ref{lm:traceetpb} is quite straightforward, see \cite[4.2.6]{SGA5}.

\section{Local terms in motivic homotopy}

\subsection{Quadratic local terms}
\label{sec:A1inv}
\subsubsection{}
As above, let $c=(c_1,c_2):C\to X\times_SX$ be a morphism, and consider the following Cartesian diagram
\begin{align}
\label{eq:fixcdiag}
\begin{split}
  \xymatrix@=10pt{
    Fix(c) \ar[r]^-{c'} \ar[d]_-{\delta'} & X \ar[d]^-{\delta_{X/S}}\\
    C \ar[r]^-{c} & X\times_SX.
  }
\end{split}
\end{align}
Let $K\in\mathcal{SH}_c(X)$ be USLA over $S$, and let $v$ be a virtual vector bundle on $C$. For a correspondence $u:c_1^*K\to c_2^!K\otimes Th(v)$ over $c$, we have defined its trace $Tr(u/S)\in H(Fix(c)/S,-v_{|Fix(c)})$ in~\eqref{eq:traceu}.

\begin{definition}
\label{def:loctrace}
For $\beta$ an open subset of $Fix(c)$, we define the element
\begin{align}
Tr_\beta(u/S)\in H(\beta/S,-v_{|\beta})
\end{align}
as the restriction of $Tr(u/S)$ to $\beta$, namely the image of $Tr(u/S)$ under the pullback map 
\begin{align}
H(Fix(c)/S,-v_{|Fix(c)})\xrightarrow{\eqref{eq:refgysin}} H(\beta/S,-v_{|\beta}). 
\end{align}
\end{definition}

\subsubsection{}
We now relate this construction with some quadratic invariants.
By a famous result of Morel, the endomorphism ring of the unit object $\mathbbold{1}_k$ is isomorphic to  the Grothendieck-Witt ring of non-degenerate symmetric bilinear forms over $k$ (or equivalently quadratic forms in characteristic different from $2$):
\begin{align}
End_{\mathcal{SH}(k)}(\mathbbold{1}_k)\simeq GW(k)
\end{align}
(see \cite{Mor}, \cite[Thm. 10.12]{BH}). This map has a twisted variant (\cite[\S5]{Mor}, \cite[\S3]{Lev}): if $v$ is a virtual vector bundle over $k$ of virtual rank $0$, then there is a canonical isomorphism
\begin{align}
H(k/k,v)\simeq GW(k,\operatorname{det}(v))
\end{align}
where the right-hand side is the twisted Grothendieck-Witt group.

\subsubsection{}
\label{num:relori}
In the setting of Definition~\ref{def:loctrace}, consider the following conditions:
\begin{enumerate}
\item
The base field $k$ has characteristic different from $2$;
\item
The scheme $S$ is smooth over $k$, and the dimension of $S$ agrees with the virtual rank of $v$;
\item
The scheme $\beta$ is proper over the base field $k$; 
\item
There exist a line bundle $w$ on $\beta$, a line bundle $v_0$ on $k$ and an isomorphism of line bundles 
\begin{align}
v_{0|\beta}
\simeq 
w^{\otimes2}\otimes\operatorname{det}((T_{S/k})_{|\beta}-v_{|\beta}).
\end{align}
\end{enumerate}
The last condition is related to \emph{relative orientations} in \cite[Def. 1.5]{BW}. Note that a case of special interest is the case where $\beta$ isomorphic the base field $k$, i.e. $\beta$ is an isolated fixed point of $c$.

Under the conditions above, and with the notations of~\eqref{eq:bivsp}, we have a map
\begin{align}
\label{eq:KOpf}
\begin{split}
H(\beta/S,-v_{|\beta})
&\simeq
H(\beta/k,(T_{S/k})_{|\beta}-v_{|\beta})
\xrightarrow{\eqref{eq:bivsec}}
KO(\beta/k,\operatorname{det}((T_{S/k})_{|\beta}-v_{|\beta})\\
&\simeq
KO(\beta/k,v_{0|\beta})
\to
KO(k/k,v_0)
\simeq
GW(k,v_0)
\end{split}
\end{align}
where $KO\in \mathcal{SH}(k)$ is the Hermitian $K$-theory spectrum which is an $SL^c$-oriented spectrum (see \cite[Appendix A]{BW}).
\begin{definition}
Under the conditions of~\ref{num:relori}, we define the \textbf{(quadratic) local term}
\begin{align}
LT_\beta(u/k)\in GW(k,v_0)
\end{align}
as the image of $Tr_\beta(u/S)$ by the map~\eqref{eq:KOpf}.
\end{definition}

\subsubsection{}
The local terms are the local contributions to the trace map, which are quadratic refinements of the classical local terms. In what follows, we show that some invariants in $\mathbb{A}^1$-enumerative geometry, such as the \emph{local $\mathbb{A}^1$-Brouwer degree} or the \emph{Euler class with support}, can be interpreted as local terms. 

\subsubsection{}
If $f:X\to S$ is a smooth morphism, then we define $M_S(X)=f_!f^!\mathbbold{1}_S$. If $U$ is an open subscheme of $X$, we define $M_S(X/U)$ to be the cofiber of $M_S(U)\to M_S(X)$ (\cite[2.3.14]{CD}). Explicitly, if $i:Z\to X$ is the immersion of the reduced closed complement of $U$, then $M_S(X/U)=f_!i_!i^*f^!\mathbbold{1}_S$.

\subsubsection{}
Let $q:X\to S$ be a smooth morphism and let $s:S\to X$ be a section of $q$. Then by \cite[Prop. 17.2.5]{EGA4} there is a canonical isomorphism of virtual vector bundles
\begin{align}
\label{eq:TNZ}
i^*\tau_{q}\simeq -\tau_s.
\end{align}
Let $c_1:C\to X$ be a morphism of smooth $S$-schemes, and consider the commutative diagram
\begin{align}
\begin{split}
  \xymatrix@=10pt{
    C_s\ar[r]^-{c_s} \ar[d]_-{s_C} \ar@{}[rd]|{\Delta} & S \ar[d]^-{s} \ar@{=}[rd] & \\
    C \ar[r]^-{c_1} \ar@/_.7pc/[rr]_-{p} & X \ar[r]^-{q} & S
  }
\end{split}
\end{align}
where the square $\Delta$ is Cartesian.
Then there is a canonical map
\begin{align}
\label{eq:map_mot_S}
M_S(C/C-C_s)\to M_S(X/X-S)
\end{align}
given by the composition
\begin{align}
\label{eq:canmapMS}
p_!s_{C!}s_C^*p^!\mathbbold{1}_S
=
c_{s!}s_C^*c_1^!q^!\mathbbold{1}_S
\xrightarrow{Ex(\Delta^{*!})}
c_{s!}c_s^!s^*q^!\mathbbold{1}_S
\to
s^*q^!\mathbbold{1}_S
=
q_!s_!s^*q^!\mathbbold{1}_S.
\end{align}
By relative purity (\cite[Cor. 2.4.37]{CD}) and the projection formula, the map~\eqref{eq:canmapMS} can be rewritten as
\begin{align}
\label{eq:purmapMS}
\begin{split}
c_{s!}Th(\tau_{c_1|C_s})\otimes Th(T_{q|S})
&\simeq
c_{s!}(Th(\tau_{c_1|C_s})\otimes Th(T_{q|C_s}))
\simeq
c_{s!}Th(T_{p|C_s})\\
&\to
s^*Th(T_q)
\simeq
Th(T_{q|S}).
\end{split}
\end{align}
Desuspending both sides by $Th(T_{q|S})$ and using adjunction, the map~\eqref{eq:purmapMS} can be rewritten as a map
\begin{align}
\label{eq:map_r2id}
\xi(\Delta):Th(\tau_{c_1|C_s})\to c_s^!\mathbbold{1}_{S}
\end{align}
which we view as a class $\xi(\Delta)\in H(C_s/S,\tau_{c_1|C_s})$.
\begin{lemma}
\label{prop:reffdl}
Under the assumptions above, the class $\xi(\Delta)$ agrees with the refined fundamental class $\Delta^*\eta_{c_1}$ in~\ref{num:DJK425}.
\end{lemma}
\proof
This follows from the construction of the map $\Delta^*$ in~\eqref{eq:verBC}, and the description of the fundamental class of the morphism $c_1$ between smooth $S$-schemes in terms of relative purity (see \cite[Thm 3.3.2]{DJK}).
\endproof

\subsubsection{}
Now consider the correspondence 
\begin{align}
c:C\xrightarrow{(c_1, s\circ p)}X\times_SX
\end{align}
with the second component given by $c_2=s\circ p:C\to X$, and we have canonically $Fix(c)=C_s$. By~\eqref{eq:TNZ}, we also have an isomorphism
\begin{align}
\label{eq:tau12}
\tau_{c_2}\simeq\tau_{c_1}+c_1^*\tau_{s\circ q}\overset{\eqref{eq:TNZ}}{\simeq}\tau_{c_1}.
\end{align}
Consider the twisted cohomological correspondence over $c$
\begin{align}
\label{eq:correfdl}
u:
c_1^*\mathbbold{1}_X
=
\mathbbold{1}_C
\overset{\eta_{c_2}}{\simeq}
c_2^!\mathbbold{1}_X\otimes Th(-\tau_{c_2})
\overset{\eqref{eq:tau12}}{\simeq}
c_2^!\mathbbold{1}_X\otimes Th(-\tau_{c_1})
\end{align}
induced by the fundamental class of $c_2$. It has a trace $Tr(u/S)\in H(C_s/S,\tau_{c_1|C_s})$ as defined in~\eqref{eq:ustrace}.

\begin{proposition}
There is a canonical homotopy between the trace $Tr(u/S)$ and the map $\xi(\Delta)$ in~\eqref{eq:map_r2id}.
\end{proposition}
\proof
Denote by $\delta=\delta_{X/S}:X\to X\times_SX$.
The trace of the map~\eqref{eq:correfdl} is given by the composition
\begin{align}
\begin{split}
Th(\tau_{c_2|C_s})
&\to
c_s^!s^!s_!c_{s!}Th(\tau_{c_2|C_s})
\simeq
c_s^!s^!\delta^*c_!Th(\tau_{c_2})
\xrightarrow{u}
c_s^!s^!\delta^*c_!c^!p_2^!\mathbbold{1}_X\\
&\to
c_s^!s^!\delta^*p_2^!\mathbbold{1}_X
\simeq
c_s^!s^!q^!\mathbbold{1}_S
=
c_s^!\mathbbold{1}_S.
\end{split}
\end{align}
The result then follows from the associativity of fundamental classes (\cite[Thm. 3.3.2]{DJK}), and the fact that the fundamental class of $s\circ q$ is invertible.
\endproof

\begin{example}[Local $\mathbb{A}^1$-Brouwer degree]
We assume $S=k$, $C=X=\mathbb{A}^n_k$, $s:k\to\mathbb{A}^n_k$ is the zero section. Let $c_1=f:\mathbb{A}^n_k\to \mathbb{A}^n_k$ be a $k$-morphism. The diagram~\eqref{eq:fixcdiag} becomes
\begin{align}
\begin{split}
  \xymatrix@=10pt{
    f^{-1}(0) \ar[r]^-{} \ar[d]_-{} & \mathbb{A}^n_k \ar[d]^-{\delta}\\
    \mathbb{A}^n_k \ar[r]^-{(f,0)} & \mathbb{A}^n_k\times_k\mathbb{A}^n_k.
  }
\end{split}
\end{align}

Assume that $x$ is an isolated zero of $f$. There is a canonical trivialization of $\tau_{c_1}$, so we have the local term 
\begin{align}
LT_{\{x\}}(u/k)\in GW(k).
\end{align}
This element recovers the \emph{local $\mathbb{A}^1$-Brouwer degree} (\cite[Def. 11]{KW}, \cite[Def. 7.1]{BW}), which we can see from the description in~\eqref{eq:map_mot_S}. 
There has been an extensive study of this element in the literature, see \cite{KW}, \cite{BBMMO}, \cite{BW} and \cite{BMP}.
\end{example}

\begin{example}[Euler class with support]
We assume $C=S$ is a smooth $k$-scheme, $X$ is a vector bundle over $S$, $s_0:S\to X$ is the zero section and $c_1=\sigma:S\to X$ is another section. The diagram~\eqref{eq:fixcdiag} becomes
\begin{align}
\begin{split}
  \xymatrix@=10pt{
    Z(\sigma) \ar[r]^-{} \ar[d]_-{} & X \ar[d]^-{\delta}\\
    S \ar[r]^-{(\sigma,s_0)} & X\times_SX
  }
\end{split}
\end{align}
where $Z(\sigma)=C_s$ is the zero locus of the section $\sigma$, and $\xi(\Delta)$ is the Euler class with support (\cite[3.2.10]{DJK}, \cite[Def. 5.13]{BW}).

Assume that $x$ is an isolated zero of the section $\sigma$, and that the dimension of $S$ agrees with the rank of the vector bundle $X$. Choose $v_0=\operatorname{det}(T_{S/k}-X)_{|x}$ and $w=0$ in~\ref{num:relori}. Then we have the local term
\begin{align}
LT_{\{x\}}(u/k)\in GW(k,\operatorname{det}(T_{S/k}-X)_{|x}).
\end{align}
This local contribution to the Euler class with support has also been studied in the literature, see \cite[\S5]{Lev} and \cite[Thm. 7.6]{BW}.
\end{example}

\subsection{Contracting correspondences and local terms}

\subsubsection{}
In this section let $k$ be a field of exponential characteristic $p$ and let $X$ be a $k$-scheme. 

\begin{definition}
\label{def:invind}
Let $c:C\to X\times_kX$ be a correspondence. A closed subset $Z\subset X$ is \textbf{$c$-invariant} if $c_1(c_2^{-1}(Z))\subset Z$. 

If $Z$ is a closed subscheme of $X$ defined by an ideal sheaf $\mathcal{I}$, we say that $c$ \textbf{stabilizes $Z$} if $c_1(c_2^{-1}(Z))\subset X$ is scheme-theoretically contained in $Z$, i.e. if $c_1^\#(\mathcal{I})\subset c_2^\#(\mathcal{I})\cdot\mathcal{O}_C$. We say that $c$ is \textbf{contracting near $Z$} if c stabilizes $Z$ and there exists $n$ such that $c_1^\#(\mathcal{I})^n\subset c_2^\#(\mathcal{I})^{n+1}\cdot\mathcal{O}_C$.

In particular, if $c$ stabilizes $Z$ then $Z$ is $c$-invariant.
\end{definition}

\subsubsection{}
The main reason why we work over a field instead of a general base lies in the following proposition:
\begin{proposition}
\label{prop:tr=0}
Assume that $Fix(c)$ is connected, and let $v$ be a virtual vector bundle over $C$.
Let $Z$ be a closed subscheme of $X$ such that $c$ is contracting near $Z$ and $c'(Fix(c))\cap Z$ is non-empty.
Let $K\in\mathcal{SH}_c(X)[1/p]$ be together with a null homotopy $K_{|Z}\simeq0$, which satisfies the condition~\ref{num:Kcondfin} below.
Then for correspondence $u:c_1^*K\to c_2^!K\otimes Th(v)$ over $c$, there is a null homotopy $Tr(u/k)\simeq0$.
\end{proposition}
The proof of Proposition~\ref{prop:tr=0} uses Ayoub's theory of motivic nearby cycles, which will be given later in~\ref{sec:tr=0}. The condition~\ref{num:Kcondfin} is satisfied for example if $k$ has characteristic $0$, or if we work with rational coefficients in $\mathcal{SH}_c(-,\mathbb{Q})$, see  Lemma~\ref{lm:condsat} below.

\subsubsection{}
\label{num:loccorr}
If $i:Z\to X$ is the inclusion of a $c$-invariant closed subscheme, let $c_{|Z}:c_2^{-1}(Z)_{red}\to Z\times_kZ$ be the restriction of $c$. Then the commutative square
\begin{align}
\label{eq:diagdeltaZ}
\begin{split}
  \xymatrix@=10pt{
    c_2^{-1}(Z)_{red} \ar[r]^-{c_{|Z}} \ar[d]_-{} \ar@{}[rd]|{\Delta_Z} & Z\times_kZ \ar[d]^-{i\times_ki}\\
    C \ar[r]_-{c} & X\times_kX
  }
\end{split}
\end{align}
satisfies the first condition of pullback in~\ref{num:corrpb}. It follows that, for every $K\in\mathcal{SH}(X)$, every virtual vector bundle $v$ over $C$ and every correspondence $u:c_1^*K\to c_2^!K\otimes Th(v)$ over $c$, the pullback by $\Delta_Z$ defines a map 
\begin{align}
u_{|\Delta_Z}:c_{|Z1}^*K_{|Z}\to c_{|Z2}^!K_{|Z}\otimes Th(v_{|c_2^{-1}(Z)_{red}}).
\end{align}
In addition, the square~\eqref{eq:diagdeltaZ} satisfies the second condition of push-forward in~\ref{num:corrpf}, so the push-forward of $u_{|\Delta_Z}$ defines a map
\begin{align}
\label{eq:deltaZ!}
\Delta_{Z!}u_{|\Delta_Z}:c_1^*i_!K_{|Z}\to c_2^!i_!K_{|Z}\otimes Th(v).
\end{align}

\subsubsection{}
\label{num:loccorr2}
On the other hand, let $j:U\to X$ be the open complement of $i$, then we have $c_1^{-1}(U)\subset c_2^{-1}(U)$. Let $c_{|U}:c_1^{-1}(U)\to U\times_kU$ be the restriction of $c$. Then the commutative square
\begin{align}
\label{eq:diagdeltaU}
\begin{split}
  \xymatrix@=10pt{
    c_1^{-1}(U) \ar[r]^-{c_{|U}} \ar[d]_-{} \ar@{}[rd]|{\Delta_U} & U\times_kU \ar[d]^-{j\times_kj}\\
    C \ar[r]_-{c} & X\times_kX
  }
\end{split}
\end{align}
satisfies the second condition of pullback in~\ref{num:corrpb}. It follows that, for every $K\in\mathcal{SH}(X)$, every virtual vector bundle $v$ over $C$ and every correspondence $u:c_1^*K\to c_2^!K\otimes Th(v)$ over $c$, the pullback by $\Delta_U$ defines a map 
\begin{align}
u_{|\Delta_U}:c_{|U1}^*K_{|U}\to c_{|U2}^!K_{|U}\otimes Th(v_{|c_1^{-1}(U)}). 
\end{align}
In addition, the square~\eqref{eq:diagdeltaZ} satisfies the first condition of push-forward in~\ref{num:corrpf}, so the push-forward of $u_{|\Delta_U}$ defines a map
\begin{align}
\label{eq:deltaU!}
\Delta_{U!}u_{|\Delta_U}:c_1^*j_!K_{|U}\to c_2^!j_!K_{|U}\otimes Th(v).
\end{align}
The following additivity result is the analogue of \cite[Prop. 1.5.10]{Var}, which we can prove using higher category theory thanks to Lemma~\ref{lm:additivity}:
\begin{lemma}
\label{lm:Var1510}
Let $K\in\mathcal{SH}_c(X)$ and let $u:c_1^*K\to c_2^!K\otimes Th(v)$ be a correspondence over $c$. Then the traces of the maps~\eqref{eq:deltaZ!} and~\eqref{eq:deltaU!} satisfy
\begin{align}
Tr(u/k)=Tr(\Delta_{Z!}u_{|\Delta_Z}/k)+Tr(\Delta_{U!}u_{|\Delta_U}/k).
\end{align}
\end{lemma}
\proof
The localization triangle gives rise to a canonical cofiber sequence in $\mathcal{SH}_c(X)$ of USLA objects over $S$
\begin{align}
j_!K_{|U}\to K\to i_!K_{|Z}
\end{align}
and the construction in~\ref{num:loccorr} and~\ref{num:loccorr2} gives a morphism of cofiber sequences
\begin{align}
\begin{split}
  \xymatrix@=11pt{
    c_1^*j_!K_{|U} \ar[r]^-{} \ar[d]_-{\Delta_{U!}u_{|\Delta_U}} & c_1^*K \ar[r]^-{} \ar[d]^-{u} & c_1^*i_!K_{|Z} \ar[d]^-{\Delta_{Z!}u_{|\Delta_Z}}\\ 
    c_2^!j_!K_{|U}\otimes Th(v) \ar[r]^-{} & c_2^!K\otimes Th(v) \ar[r]^-{} & c_2^!i_!K_{|Z}\otimes Th(v)
  }
\end{split}
\end{align}
in the $\infty$-categorical sense. We conclude by applying Lemma~\ref{lm:additivity}.
\endproof
\subsubsection{}
The following lemma is proved in \cite[Thm. 2.1.3 (a)]{Var}:
\begin{lemma}
\label{lm:icrediso}
Assume that $c$ is contracting near $Z$, $Fix(c)$ is connected, and $c'(Fix(c))\cap Z$ is non-empty.



Then the canonical closed immersion $i_c:Fix(c_{|Z})\to Fix(c)$ induces an isomorphism $(i_c)_{red}:(Fix(c_{|Z}))_{red}\simeq (Fix(c))_{red}$. 

In particular, the push-forward map $i_{c*}:H(Fix(c_{|Z})/k,-v_{|Fix(c_{|Z})})\to H(Fix(c)/k,-v_{|Fix(c)})$ is an isomorphism.
\end{lemma}

\subsubsection{}
We now deal with the general case of not necessarily $c$-invariant subschemes. Let $v$ be a virtual vector bundle over $C$.
\begin{definition}
\label{def:contrnbh}
\begin{enumerate}
\item
Let $Z$ be a closed subscheme of $X$. Let $W$ be the complement of the closure of $c_2^{-1}(Z)-c_1^{-1}(Z)$ in $C$. Then $W$ is the largest subset of $C$ such that $Z$ is $c_{|W}$-invariant (\cite[Lemma 1.5.3]{Var}). We still denote by $c_{|W}:W\to X\times_kX$ be the restriction, and let
\begin{align}
c_{|Z}=(c_{|W})_{|Z}:c_{|W2}^{-1}(Z)_{red}\to Z\times_kZ
\end{align}
be defined as in~\ref{num:loccorr}.

For $K\in\mathcal{SH}(X)$, $v$ a virtual vector bundle over $C$ and $u:c_1^*K\to c_2^!K\otimes Th(v)$ a correspondence over $c$, let $u_{|\Delta_W}:c_{|W1}^*K\to c_{|W2}^!K\otimes Th(v_{|W})$ be the pullback of $u$, and let
\begin{align}
u_{Z}=(u_{|\Delta_W})_{|\Delta_Z}:c_{|Z1}^*K_{|Z}\to c_{|Z2}^!K_{|Z}\otimes Th(v_{|c_{|W2}^{-1}(Z)_{red}})
\end{align}
be defined as in~\ref{num:loccorr}, considered as a correspondence over $c_{|Z}$.

\item
We say that $c$ is \textbf{contracting near $Z$ in a neighborhood of fixed points} if there is an open neighborhood $W\subset C$ of $Fix(c)$ such that $c_{|W}$ is contracting near $Z$ in the sense of Definition~\ref{def:invind}.

\end{enumerate}
\end{definition}

\subsubsection{}
\label{num:contrnbh}
Let $c:C\to X\times_kX$ be a correspondence contracting near a closed subscheme $Z$ in a neighborhood of fixed points. Let $\beta$ be an open connected subset of $Fix(c)$ such that $c'(\beta)\cap Z$ is non-empty. Then by Lemma~\ref{lm:icrediso}, there is a unique open connected subscheme $\beta'$ of $Fix(c_{|Z})$ such that $\beta=i_c(\beta')$.

\subsubsection{}
In the following theorem, similar to the case of Proposition~\ref{prop:tr=0}, we need some condition related to the motivic nearby cycle which will be stated in~\ref{num:Kcondfin} below. Such a condition is satisfied if the base field $k$ has characteristic $0$, or if we work with rational coefficients, see  Lemma~\ref{lm:condsat} below.
\begin{theorem}
\label{th:LTnaive}
Assume the setting of~\ref{num:contrnbh}.
Let $K\in\mathcal{SH}_c(X)[1/p]$ which satisfies the condition~\ref{num:Kcondfin} below.
Let $v$ be a virtual vector bundle over $C$ and let $u:c_1^*K\to c_2^!K\otimes Th(v)$ be a correspondence over $c$. Then we have 
\begin{align}
Tr_\beta(u/k)=i_{c*}Tr_{\beta'}(u_{Z}/k)\in H(Fix(c)/k,-v_{|\beta}). 
\end{align}
In particular, if we further assume that the conditions of~\ref{num:relori} are satisfied. Then we have
\begin{align}
LT_\beta(u/k)=LT_{\beta'}(u_{Z}/k)\in GW(k,v_0).
\end{align}
\end{theorem}
\proof
Let $W\subset C$ be an open neighborhood of $Fix(c)$ such that $c_{|W}$ is contracting near $Z$, then $Fix(c_{|W})=Fix(c)$. Therefore by replacing $c$ by $c_{|W}$ and $u$ by $u_{|\Delta_W}$, we may assume that $c$ is contracting near $Z$. Also by replacing $C$ by the open subscheme $C-(Fix(c)-\beta)$, we may assume that $\beta=Fix(c)$. By Lemma~\ref{lm:Var1510} and Lemma~\ref{lm:LZ221} we have
\begin{align}
Tr_\beta(u)
=
Tr_\beta(\Delta_{Z!}u_{|\Delta_Z})+Tr_\beta(\Delta_{U!}u_{|\Delta_U})
=
i_{c*}Tr_{\beta'}(u_{Z})+Tr_\beta(\Delta_{U!}u_{|\Delta_U}).
\end{align}
By Proposition~\ref{prop:tr=0} we have $Tr_\beta(\Delta_{U!}u_{|\Delta_U})\simeq0$, which finishes the proof.
\endproof

\begin{remark}
The element $LT_{\beta'}(u_{Z}/k)$ is called the \textbf{naive local term}. For example (see \cite[Example 1.5.7]{Var}), if the morphism $c_2$ is quasi-finite, then for each closed point $x$ of $X$ the set $Fix(c_{|x})=c_1^{-1}(x)\cap c_2^{-1}(x)$ is finite. Each point $y\in Fix(c_{|x})$ determines an endomorphism $u_y:K_{|x}\to K_{|x}$, and we have $LT_{y}(u_{x}/k)=Tr(u_y/k)$ is the usual (categorical) trace.

\end{remark}

\section{Specialization and deformation of correspondences}
\label{sec:spdef}
This section is devoted to the proof of Proposition~\ref{prop:tr=0} as promised.

\subsection{The motivic nearby cycle functor}
\subsubsection{}
Let $k$ be a field of exponential characteristic $p$ and let $\sigma=\operatorname{Spec}(k)\to S=\mathbb{A}^1_k$ be the inclusion of the zero section, with open complement $\eta=\mathbb{G}_{m,k}$. For every morphism $f:X\to S$, the \textbf{(tame) motivic nearby cycle functor} is a functor 
\begin{align}
\Psi_X:\mathcal{SH}(X_\eta)\to\mathcal{SH}(X_\sigma).
\end{align}
Ayoub defined this functor in \cite[Def. 3.5.6]{Ayo} in the language of $1$-categories and derivators, which can be enhanced into a functor of $\infty$-categories, see~\cite{JY2}. 

\subsubsection{}
\label{num:psiprop}
Ayoub has shown that the functor $\Psi$ satisfies some good cohomological properties (\cite[\S3.5]{Ayo}): 
\begin{enumerate}
\item
(\cite[Prop. 3.2.17]{Ayo}) 
The functor $\Psi_X$ is pseudo-monoidal, that is, there is a binatural transformation
\begin{align}
\label{eq:pmonpsi}
\Psi_X(K)\otimes\Psi_X(L)\to\Psi_X(K\otimes L).
\end{align}
Consequently, for two $S$-schemes $X$ and $Y$, there is a binatural tranformation
\begin{align}
\label{eq:kunpsi}
\Psi_X(K)\boxtimes_\sigma\Psi_Y(L)\to\Psi_{X\times_SY}(K\boxtimes_\eta L).
\end{align}
\item
(\cite[Prop. 3.1.9, 3.2.9]{Ayo})
For every morphism $g:Y\to X$, there are natural transformations
\begin{align}
\label{eq:psiupper*}
g_\sigma^*\circ \Psi_X\to \Psi_Y\circ g_\eta^*
\end{align}
\begin{align}
\label{eq:psiupper!}
\Psi_Y\circ g_\eta^!\to g_\sigma^!\circ \Psi_X
\end{align}
which are invertible if $g$ is smooth. By adjunction, the map~\eqref{eq:psiupper*} induces a natural transformation
\begin{align}
\label{eq:psilower*}
\Psi_X\circ g_{\eta*}\to g_{\sigma*}\circ \Psi_Y
\end{align}
which is invertible when $g$ is proper.
\item 
(\cite[Def. 3.2.11]{Ayo}) 
There is a canonical natural transformation 
\begin{align}
\label{eq:chipsi}
i_X^*j_{X*}\to \Psi_X
\end{align}
where $i_X:X_\sigma\to X$ and $j_X:X_\eta\to X$ are the immersions, which arises as a natural transformation of specializations systems.
\item
(\cite[Lemme 3.5.10]{Ayo}) 
There is a canonical natural transformation of functors 
\begin{align}
\label{eq:Ayo3510}
1
\to
i^*j_*q^*
\to
\Psi_Sq^*
\end{align}
which is an isomorphism, where $q:\eta\to\operatorname{Spec}(k)$ is the projection.

\item
(\cite[Prop. 3.1.7]{Ayo}) 
The functor $\Psi$ is compatible with the formation of Thom spaces of vector bundles: if $X$ is an $S$-scheme and $v$ is a virtual vector bundle over $X$, then there is a canonical isomorphism 
\begin{align}
\label{eq:spthom}
Th(v_{|X_\sigma})\otimes\Psi_X(-)\simeq\Psi_X(-\otimes Th(v_{|X_\eta})).
\end{align}
In particular, the functor $\Psi$  commute with Tate twists.
\end{enumerate}

\begin{remark}
When $k$ has characteristic $0$, the K\"unneth-type map~\eqref{eq:kunpsi} is invertible (\cite[Cor. 3.5.18]{Ayo}). This property no longer holds in positive characteristics, already in \'etale cohomology. We refer to \cite{Ill} for a detailed discussion about K\"unneth formulas for nearby cycles. In \cite[Thm. 10.19]{Ayo2} Ayoub proves a similar K\"unneth formula for the \emph{total nearby cycle functor} in \'etale motives. In what follows, we only need the existence of the map~\eqref{eq:kunpsi} and not the K\"unneth property.
\end{remark}

\subsubsection{}
\label{num:propnearby}
We now discuss some more technical properties. 
By \cite[Thm. 3.5.14, 3.5.20]{Ayo}, we know that there is a natural transformation
\begin{align}
\label{eq:Psidual}
\Psi_X\mathbb{D}_{X_\eta/\eta}\to \mathbb{D}_{X_\sigma/\sigma}\Psi_X.
\end{align}
Now let $\mathfrak{P}$ be a collection of primes and let $K\in\mathcal{SH}_c(X_\eta)[\mathfrak{P}^{-1}]$. Consider the following properties:
\begin{enumerate}
\item
The object $\Psi_X(K)\in\mathcal{SH}(X_\sigma)[\mathfrak{P}^{-1}]$ is constructible;
\item
The map~\eqref{eq:Psidual} is invertible when applied to $K$.
\end{enumerate}
We now discuss cases where the conditions of~\ref{num:propnearby} are fulfilled.

\begin{definition}[\textrm{\cite[Def. 3.3.33]{Ayo}}]
\label{def:Ayo3333}
Let $f:X\to S$ be a morphism of finite type and let $x$ be a point of $X_\sigma$. 
We say that $f$ has \textbf{globally semi-stable reduction at $x$} (with $n$ branches of type $\underline{a}=(a_1,\cdots,a_n)\in\mathbb{N}_{>0}^n$) if there exists a Zariski neighborhood $U$ of $x$ in $X$ such that
\begin{enumerate}
\item $U$ is regular; 

\item There exists $n+1$ global sections $t_1,\cdots,t_n,u$ of $U$ such that
\begin{enumerate}
\item $u$ is an invertible global section and $\pi=u\cdot t_1^{a_1}\cdots t_n^{a_n}$, where $\pi:S\to\mathbb{A}^1_k$ is the identity morphism, considered as a global section of $X$;
\item For all $1\leqslant i\leqslant n$, the closed subscheme $D_i$ of $U$ defined by $t_i=0$ is regular and contains $x$;
\item The union $\cup_{i=1}^nD_i$ is a strict normal-crossing divisor in $U$.
\end{enumerate}
\end{enumerate}
We say that $f$ has \textbf{globally good semi-stable reduction at $x$} if in addition the following condition is satisfied:
\begin{enumerate}
\setcounter{enumi}{2}
\item
\label{num:cond3}
Let $p$ be the exponential characteristic of $k$.
Let $m_0\in\mathbb{N}\cup\{+\infty\}$ be the upper bound of all non-negative integers $m$ such that there exists a global section $v$ of $U$ with $v^{p^m}=u$. If $m_0=+\infty$, we set $v_0=u$, otherwise let $v_0$ be a global section of $U$ such that $v_0^{p^{m_0}}=u$. Let $l_0$ be the integer given by the upper bound
\begin{align}
l_0=\operatorname{sup}\left(\{l\in\mathbb{N}|\forall i\in\{1,\cdots,n\},p^{m_0+l}\textrm{ divides }a_i\}\cup\{0\}\right).
\end{align}
Let $O$ be the closed subscheme of $U$ defined by the equations $t_1=\cdots=t_n=0$. Then for any $l\in\{0,\cdots,l_0\}$, the scheme $O[v_0^{\frac{1}{p^l}}]$ is smooth over $k$.
\end{enumerate}

We say that $f$ or $X$ has \textbf{globally semi-stable reduction} (respectively, \textbf{globally good semi-stable reduction}) if there exists an integer $N$ such that for any point $x$ of $X_\sigma$, $f$ has globally semi-stable reduction at $x$ with $n_x$ branches with $n_x\leqslant N$.

\end{definition}

\subsubsection{}
According to our conventions, if the scheme $X_\sigma$ is empty then $X$ has globally good semi-stable reduction. If $k$ has characteristic $0$, then globally semi-stable reduction is equivalent to globally good semi-stable reduction, since the condition~\eqref{num:cond3} above is trivial.

\begin{corollary}
\label{cor:genss}
Let $f:X\to S$ be a morphism of finite type and let $\mathfrak{P}$ be a collection of primes containing $p$.
Then $\mathcal{SH}_c(X_\eta)[\mathfrak{P}^{-1}]$ is generated as a thick subcategory by elements of the form $g_{\eta*}\mathbbold{1}_{X'_\eta}(n)$, where 
\begin{itemize}
\item $g:X'\to X$ is a projective morphism;
\item $X'$ is smooth over a finite purely inseparable extension of $k$;
\item $f\circ g$ has globally semi-stable reduction;
\item $n$ is an integer.
\end{itemize}

\end{corollary}
\proof
As in \cite[Thm. 3.3.46]{Ayo}, applying Lemma~\ref{lm:CD7.2} to $X$ with the closed subset $X_\sigma$ and the motivic $\infty$-category 
\begin{align}
\begin{split}
Sch/S&\to Pr^L_{st}\\ 
X&\mapsto \mathcal{SH}(X_\eta)[\mathfrak{P}^{-1}]
\end{split}
\end{align}
we know that $\mathcal{SH}_c(X_\eta)[\mathfrak{P}^{-1}]$ is generated as a thick subcategory by elements of the form $g_{\eta*}\mathbbold{1}_{X'_\eta}(n)$, where $g:X'\to X$ is a projective morphism, $X'$ is smooth over a finite purely inseparable extension of $k$, $n$ is an integer, and $(X'_{\sigma})_{red}$ is either empty, or a strict normal-crossing divisor on $X'$. It remains to show that if $(X'_{\sigma})_{red}$ is a strict normal-crossing divisor on $X'$, then $X'$ has globally semi-stable reduction in the sense of Definition~\ref{def:Ayo3333}. 

We use the arguments in the proof of \cite[Prop. 2.5]{AIS}. Let $x\in X'_{\sigma}$ be a point. Then there exists an affine neighborhood $U$ of $x$ in $X'$ such that each of the irreducible components $D_1,\cdots,D_n$ of the divisor $(U_\sigma)_{red}=(X_\sigma)_{red}\cap U$ contains $x$, and is defined by a single equation. Let $t_i\in\mathcal{O}(U)$ be an equation defining $D_i$, and let $a_i$ be the multiplicity of $D_i$ in $U_\sigma$. Then $t_1^{a_1}\cdots t_n^{a_n}$ is a generator of the ideal defining $U_\sigma$ in $U$. The latter ideal is also generated by the uniformizer $\pi$, which implies that there exists an invertible section $u\in\mathcal{O}^*(U)$ such that $\pi=u\cdot t_1^{a_1}\cdots t_n^{a_n}$, and therefore $X'$ has globally semi-stable reduction.
\endproof

\begin{definition}
Let $\mathfrak{P}$ be a collection of primes containing $p$.
We let $\mathcal{SH}^{good}_c(X_\eta)[\mathfrak{P}^{-1}]$ be the thick subcategory of $\mathcal{SH}_c(X_\eta)[\mathfrak{P}^{-1}]$ generated by elements of the form $g_{\eta*}\mathbbold{1}_{X'_\eta}(n)$, where 
\begin{itemize}
\item $g:X'\to X$ is a projective morphism;
\item $X'$ is smooth over a finite purely inseparable extension of $k$;
\item $f\circ g$ has globally good semi-stable reduction;
\item $n$ is an integer.
\end{itemize}
\end{definition}

\begin{lemma}
\label{lm:condsat}
Let $\mathfrak{P}$ be a collection of primes containing $p$.
The conditions of~\ref{num:propnearby} are satisfied in the following cases:
\begin{enumerate}
\item
The object $K$ belongs to $\mathcal{SH}^{good}_c(X_\eta)[\mathfrak{P}^{-1}]$.
\item
$\mathfrak{P}$ is the collection of all primes, and 
the object $K$ belongs to $\mathcal{SH}_c(X_\eta,\mathbb{Q})$.
\end{enumerate}
In particular, if $k$ has characteristic $0$, then this property holds for any $K\in\mathcal{SH}_c(X_\eta)$.
\end{lemma}
\proof
The first case can be proved with arguments in the proof of \cite[Thm. 3.5.14, 3.5.20]{Ayo}, using the variants of Thm. 3.3.46 and Cor. 3.3.49 of loc. cit.. For the second case, we have a decomposition $\mathcal{SH}_\mathbb{Q}\simeq\mathcal{DM}_\mathbb{Q}\oplus\mathcal{SH}_{\mathbb{Q},-}$, where $\mathcal{DM}_\mathbb{Q}$ satisfies \'etale descent, and $\mathcal{SH}_{\mathbb{Q},-}$ is concentrated on the characteristic $0$ fiber (see \cite[Lemma 2.2]{DFJK}). The case of $\mathcal{DM}_\mathbb{Q}$ follows from \cite[\S10]{Ayo2}, where the condition on the cohomological dimension in 7.3 of loc. cit. can be dropped by working with locally constructible $h$-motives, see \cite{CD2} and \cite{Cis} for more details; the case of $\mathcal{SH}_{\mathbb{Q},-}$ follows from the characteristic $0$ case, which proves the result.
\endproof

\subsubsection{}
We end this section with a discussion on the specialization maps.
Consider a Cartesian diagram
\begin{align}
\begin{split}
  \xymatrix@=11pt{
    X_\sigma \ar[r]^-{i_X} \ar[d]_-{f_\sigma} & X \ar[d]^-{f} & X_\eta \ar[d]^-{f_\eta} \ar[l]_-{j_X}\\
    \sigma \ar[r]^-{i} & S & \eta.  \ar[l]_-{j}
  }
\end{split}
\end{align}
We have a canonical map
\begin{align}
\label{eq:sppsi}
\begin{split}
j_{X*}f_\eta^!\mathbbold{1}_\eta
\to
i_{X*}i_X^*j_{X*}f_\eta^!\mathbbold{1}_\eta
\xrightarrow{\eqref{eq:chipsi}}
i_{X*}\Psi_X(f_\eta^!\mathbbold{1}_\eta)
\xrightarrow{\eqref{eq:psiupper!}}
i_{X*}f_\sigma^!\Psi_S(\mathbbold{1}_\eta)
\overset{\eqref{eq:Ayo3510}}{\simeq}
i_{X*}f_\sigma^!\mathbbold{1}_\sigma
\end{split}
\end{align}
from which we deduce, for $v$ a virtual vector bundle over $X$, another specialization map
\begin{align}
\label{eq:spnearby}
H(X_\eta/\eta,-v_{|X_\eta})\to H(X_\sigma/\sigma,-v_{|X_\sigma}).
\end{align}
The following lemma is stated in \cite[Rmk. 4.5.8]{DJK} without proof, which we will prove now.
\begin{lemma}
\label{lm:speq}
There is a canonical homotopy between the maps~\eqref{eq:spnearby} and~\eqref{eq:DJK456b}.
\end{lemma}
\proof
The map~\eqref{eq:DJK456b} is induced by the natural transformation~\eqref{eq:DJK456a}, which comes from the functor $j_{X*}j_X^*$ as the cofiber of the map $i_{X*}i_X^!\to1$ and the null-homotopy of the composition~\eqref{eq:DJKspint}. Therefore it suffices to show that the following diagram commutes:
\begin{align}
\begin{split}
  \xymatrix@=11pt{
    f^!\mathbbold{1} \ar[r]^-{} \ar[d]_-{} & i_{X*}i_X^*f^!\mathbbold{1} \ar[r]^-{\eqref{eq:purtr}} \ar[d]^-{} & i_{X*}(i_X^!f^!\mathbbold{1}\otimes Th(N_\sigma S_{|X_\sigma})) \ar[r]_-{} & i_{X*}f_\sigma^!\mathbbold{1} \\
    j_{X*}j_X^*f^!\mathbbold{1} \ar[r]^-{} & i_{X*}i_X^*j_{X*}f_\eta^!\mathbbold{1} \ar[r]^-{} & i_{X*}\Psi_X(f_\eta^!\mathbbold{1}) \ar[r]_-{} & i_{X*}f_\sigma^!\Psi_S(\mathbbold{1}). \ar[u]^-{\wr}
  }
\end{split}
\end{align}
The square on the left clearly commutes. For the diagram on the right, we are reduced to the following one:
\begin{align}
\label{eq:htpsp}
\begin{split}
  \xymatrix@=10pt{
     & i_X^*f^!\mathbbold{1} \ar[r]^-{\eqref{eq:purtr}} \ar[ld]^-{} \ar[rd]_-{} & i_X^!f^!\mathbbold{1}\otimes Th(N_\sigma S_{|X_\sigma}) \ar[d]_-{} \\
    i_X^*j_{X*}f_\eta^!\mathbbold{1} \ar[r]^-{} \ar[d]^-{} & f_\sigma^!i^*j_*\mathbbold{1} \ar[d]_-{} & f_\sigma^!\mathbbold{1} \ar[l]_-{} \ar[ld]_-{\sim} \\
    \Psi_X(f_\eta^!\mathbbold{1}) \ar[r]^-{} & f_\sigma^!\Psi_S(\mathbbold{1}). &
  }
\end{split}
\end{align}
For the commutativity of diagram~\eqref{eq:htpsp}: the upper square is clear; the upper right triangle follows from the compatibility of the Gysin morphism; the lower right triangle follows from \cite[Lemme 3.5.10]{Ayo}; the lower square follows from the natural transformation of specialization systems $i^*j_*\to\Psi$ in~\eqref{eq:chipsi}.
\endproof

\subsection{Specialization of correspondences}
\label{sec:spcorr}

\subsubsection{}
Let $f:X\to S$ and $c:C\to X\times_SX$ be two morphisms. 
We have the following Cartesian diagrams:
\begin{align}
\begin{split}
  \xymatrix@=10pt{
    C_\sigma \ar[r]^-{i_C} \ar[d]_-{c_\sigma} & C \ar[d]^-{c} & C_\eta \ar[l]_-{j_C} \ar[d]^-{c_\eta}\\
    X_\sigma\times_SX_\sigma \ar[r]^-{} \ar[d]_-{} & X\times_SX \ar[d]_-{} & X_\eta\times_\eta X_\eta \ar[l]_-{} \ar[d]_-{} \\
    \sigma \ar[r]^-{i} & S & \eta \ar[l]_-{j}
  }
\end{split}
\end{align}
\begin{align}
\begin{split}
  \xymatrix@=11pt{
    X_\sigma \ar[r]^-{i_X} \ar[d]_-{f_\sigma} & X \ar[d]^-{f} & X_\eta \ar[d]^-{f_\eta} \ar[l]_-{j_X}\\
    \sigma \ar[r]^-{i} & S & \eta \ar[l]_-{j}.
  }
\end{split}
\end{align}
Let $v$ be a virtual vector bundle on $C$. Then for any $K\in\mathcal{SH}(X_\eta)$ and any correspondence $u:c_{\eta1}^*K\to c_{\eta2}^!K\otimes Th(v_{|C_\eta})$, we have the following composition
\begin{align}
\label{eq:spcorr}
\begin{split}
c_{\sigma1}^*\Psi_X(K)
&\xrightarrow{\eqref{eq:psiupper*}}
\Psi_C(c_{\eta1}^*K)
\xrightarrow{u}
\Psi_C(c_{\eta2}^!K\otimes Th(v_{|C_\eta}))\\
&\simeq
\Psi_C(c_{\eta2}^!K)\otimes Th(v_{|C_\sigma})
\xrightarrow{\eqref{eq:psiupper!}}
c_{\sigma2}^!\Psi_X(K)\otimes Th(v_{|C_\sigma})
\end{split}
\end{align}
where the isomorphism follows from the compatibility of $\Psi$ with Thom spaces (\eqref{eq:spthom}). This construction gives rise to a canonical map
\begin{align}
\label{eq:spcorr}
Map(c_{\eta1}^*K, c_{\eta2}^!K\otimes Th(v_{|C_\eta}))\to Map(c_{\sigma1}^*\Psi_X(K), c_{\sigma2}^!\Psi_X(K)\otimes Th(v_{|C_\sigma}))
\end{align}
called the \textbf{specialization of correspondences}.

\subsubsection{}
There is a Cartesian diagram
\begin{align}
\begin{split}
  \xymatrix@=10pt{
    Fix(c_\sigma) \ar[r]^-{} \ar[d]_-{} & Fix(c) \ar[d]^-{} & Fix(c_\eta) \ar[d]^-{} \ar[l]_-{}\\
    \sigma \ar[r]^-{i} & S & \eta \ar[l]_-{j}
  }
\end{split}
\end{align}
so by~\eqref{eq:DJK456b} there is a specialization map
\begin{align}
\label{eq:spFix}
H(Fix(c_\eta)/\eta,-v_{|Fix(c_\eta)})\to H(Fix(c_\sigma)/\sigma,-v_{|Fix(c_\sigma)}).
\end{align}

\subsubsection{}
The following lemma states that specializations are compatible with the trace map:
\begin{lemma}
\label{lm:sptrace}
Let $K\in\mathcal{SH}_c(X_\eta)[1/p]$. Assume further that $K$ is USLA over $\eta$, and satisfies the conditions of~\ref{num:propnearby}.
Then the following diagram is commutative:
\begin{align}
\begin{split}
  \xymatrix{
    Map(c_{\eta1}^*K,c_{\eta2}^!K\otimes Th(v_{|C_\eta})) \ar[r]^-{Tr(-/\eta)} \ar[d]_-{\eqref{eq:spcorr}} & H(Fix(c_\eta)/\eta,-v_{|Fix(c_\eta)}) \ar[d]^-{\eqref{eq:spFix}}\\
    Map(c_{\sigma1}^*\Psi_X(K), c_{\sigma2}^!\Psi_X(K)\otimes Th(v_{|C_\sigma})) \ar[r]^-{Tr(-/\sigma)} & H(Fix(c_\sigma)/\sigma,-v_{|Fix(c_\sigma)}).
  }
\end{split}
\end{align}
where the lower horizontal map is well-defined since $\Psi_X(K)\in\mathcal{SH}(X_\sigma)[1/p]$ is USLA over $\sigma$, which follows from by Lemma~\ref{lm:USLA} as $\sigma$ is the spectrum of a field. 
\end{lemma}
\proof
Let $\phi:Fix(c)\to S$ be the structure morphism.
By Lemma~\ref{lm:speq}, the map~\eqref{eq:spFix} agrees with the map~\eqref{eq:spnearby}, which is induced by the map~\eqref{eq:sppsi}. With the notations from Section~\ref{sec:corrtr}, we are reduced to show the commutativity of the following diagram:
\begin{align}
\begin{split}
  \xymatrix@=10pt{
    Th(-v_{|Fix(c)}) \ar[r]^-{Tr(u/\eta)} \ar[d]_-{} & j_{Fix(c),*}\phi_\eta^!\mathbbold{1}_\eta \ar[d]^-{}\\
    i_{Fix(c),*}c'^!_\sigma\delta_{X_\sigma/\sigma}^*c_{\sigma!}Th(-v_{|C_\sigma}) \ar[d]_-{u'} & i_{Fix(c),*}i_{Fix(c)}^*j_{Fix(c),*}\phi_\eta^!\mathbbold{1}_\eta \ar[d]_-{} \\
    i_{Fix(c),*}c'^!_\sigma\delta_{X_\sigma/\sigma}^*c_{\sigma!}\underline{Hom}(c_{\sigma1}^*\Psi_X(K),c_{\sigma2}^!\Psi_X(K)) \ar[d]_-{\wr} & i_{Fix(c),*}\Psi_{Fix(c)}(\phi_\eta^!\mathbbold{1}_\eta) \ar[d]_-{} \\
    i_{Fix(c),*}c'^!_\sigma\delta_{X_\sigma/\sigma}^*c_{\sigma!}c_\sigma^!(\mathbb{D}_{X_\sigma/\sigma}(\Psi_X(K))\boxtimes_S\Psi_X(K)) \ar[d]_-{} & i_{Fix(c),*}\phi_\sigma^!\Psi_S(\mathbbold{1}_\eta) \ar[d]_-{\wr} \\
    i_{Fix(c),*}c'^!_\sigma\delta_{X_\sigma/\sigma}^*(\mathbb{D}_{X_\sigma/\sigma}(\Psi_X(K))\boxtimes_S\Psi_X(K)) \ar[d]_-{\wr} & i_{Fix(c),*}\phi_\sigma^!\mathbbold{1}_\sigma \ar[d]_-{\wr} \\
    i_{Fix(c),*}c'^!_\sigma(\mathbb{D}_{X_\sigma/\sigma}(\Psi_X(K))\otimes\Psi_X(K)) \ar[r]_-{} & i_{Fix(c),*}c'^!_\sigma\mathcal{K}_{X_\sigma/\sigma}
  }
\end{split}
\end{align}
where $u'$ is the specialization of $u$ defined in~\eqref{eq:spcorr}. By virtue of the commutative diagram
\begin{align}
\begin{split}
  \xymatrix@=10pt{
    1 \ar[r]^-{} \ar[d]_-{} & i_{Fix(c),*}i_{Fix(c)}^*j_{Fix(c),*}c'^!_\eta c'_{\eta!}j_{Fix(c)}^* \ar[d]^-{}\\
    i_{Fix(c),*}c'^!_\sigma c'_{\sigma!}i_{Fix(c)}^* \ar[d]_-{u'} & i_{Fix(c),*}\Psi_{Fix(c)}c'^!_\eta c'_{\eta!}j_{Fix(c)}^* \ar[d]_-{} \\
    i_{Fix(c),*}c'^!_\sigma i_X^*j_{X*}c'_{\eta!}j_{Fix(c)}^* \ar[r]_-{} & i_{Fix(c),*}c'^!_\sigma\Psi_{X}c'_{\eta!}j_{Fix(c)}^*
  }
\end{split}
\end{align}
we are reduced to the commutativity of the following diagram
\begin{align}
\begin{split}
  \xymatrix@=10pt{
    c'_{\sigma!}i_{Fix(c)}^*Th(-v_{|Fix(c)}) \ar[r]^-{} \ar[d]_-{\wr} & i_X^*j_{X*}c'_{\eta!}j_{Fix(c)}^*Th(-v_{|Fix(c)}) \ar[d]^-{}\\
    \delta_{X_\sigma/\sigma}^*c_{\sigma!}Th(-v_{|C_\sigma}) \ar[d]_-{u'} & \Psi_Xc'_{\eta!}j_{Fix(c)}^*Th(-v_{|Fix(c)}) \ar[d]^-{Tr(u/\eta)} \\
    \delta_{X_\sigma/\sigma}^*c_{\sigma!}\underline{Hom}(c_{\sigma1}^*\Psi_X(K),c_{\sigma2}^!\Psi_X(K)) \ar[d]_-{\wr} & \Psi_Xc'_{\eta!}\phi_\eta^!\mathbbold{1}_\eta \ar[d]_-{} \\
    \delta_{X_\sigma/\sigma}^*c_{\sigma!}c_\sigma^!(\mathbb{D}_{X_\sigma/\sigma}(\Psi_X(K))\boxtimes_S\Psi_X(K)) \ar[d]_-{} & \Psi_Xf_\eta^!\mathbbold{1}_\eta \ar[d]_-{} \\
    \delta_{X_\sigma/\sigma}^*(\mathbb{D}_{X_\sigma/\sigma}(\Psi_X(K))\boxtimes_S\Psi_X(K)) \ar[d]_-{\wr} & f_\sigma^!\Psi_S\mathbbold{1}_\eta \ar[d]_-{\wr} \\
    \mathbb{D}_{X_\sigma/\sigma}(\Psi_X(K))\otimes\Psi_X(K))\ar[r]_-{} & \mathcal{K}_{X_\sigma/\sigma}.
  }
\end{split}
\end{align}
This diagram is reduced to the following one:
\begin{center}
$$
\resizebox{1.1\textwidth}{!}{
  \xymatrix{
  & \delta^*_{X_\sigma/\sigma}c_{\sigma!}i_C^*j_{C*}Th(-v_{|C_\eta}) \ar[r]^-{} \ar[d]^-{} & \delta^*_{X_\sigma/\sigma}i_{X\times_SX}^*j_{X\times_SX,*}c_{\eta!}Th(-v_{|C_\eta}) \ar[r]^-{} \ar[d]^-{} & i_X^*j_{X*}\delta_{X_\eta/\eta}^*c_{\eta!}Th(-v_{|C_\eta}) \ar[d]^-{} \\
 \delta^*_{X_\sigma/\sigma}c_{\sigma!}Th(-v_{|C_\sigma}) \ar[ru]^-{} \ar[dd]^-{u'} & \delta^*_{X_\sigma/\sigma}c_{\sigma!}\Psi_{C}Th(-v_{|C_\eta}) \ar[r]^-{} \ar[d]^-{u} & \delta^*_{X_\sigma/\sigma}\Psi_{X\times_SX}c_{\eta!}Th(-v_{|C_\eta}) \ar[r]^-{} \ar[d]^-{u} & \Psi_X\delta_{X_\eta/\eta}^*c_{\eta!}Th(-v_{|C_\eta}) \ar[d]^-{u} \\
  & \delta^*_{X_\sigma/\sigma}c_{\sigma!}\Psi_{C}\underline{Hom}(c_{\eta1}^*K,c_{\eta2}^!K) \ar[r]^-{} \ar[d]^-{\wr} \ar[ld]^-{} & \delta^*_{X_\sigma/\sigma}\Psi_{X\times_SX}c_{U!}\underline{Hom}(c_{\eta1}^*K,c_{\eta2}^!K) \ar[r]^-{} \ar[d]^-{\wr} & \Psi_X\delta_{X_\eta/\eta}^*c_{\eta!}\underline{Hom}(c_{\eta1}^*K,c_{\eta2}^!K) \ar[d]^-{\wr} \\
 \delta^*_{X_\sigma/\sigma}c_{\sigma!}\underline{Hom}(c_{\sigma1}^*\Psi_X(K),c_{\sigma2}^!\Psi_X(K)) \ar[d]^-{\wr} & \delta^*_{X_\sigma\sigma}c_{\sigma!}\Psi_{C}c_\eta^!(\mathbb{D}_{X_\eta/\eta}(K)\boxtimes_SK) \ar[r]^-{} \ar[d]^-{} & \delta^*_{X_\sigma/\sigma}\Psi_{X\times_SX}c_{\eta!}c_\eta^!(\mathbb{D}_{X_\eta/\eta}(K)\boxtimes_SK) \ar[r]^-{} \ar[d]^-{} & \Psi_X\delta_{X_\eta/\eta}^*c_{\eta!}c_U^!(\mathbb{D}_{X_\eta/\eta}(K)\boxtimes_SK) \ar[d]^-{} \\
 \delta^*_{X_\sigma/\sigma}c_{\sigma!}c_\sigma^!(\mathbb{D}_{X_\sigma/\sigma}(\Psi_X(K))\boxtimes_S\Psi_X(K)) \ar[r]^-{} \ar[rd]^-{} & \delta^*_{X_\sigma/\sigma}c_{\sigma!}c_\sigma^!\Psi_{X\times_SX}(\mathbb{D}_{X_\eta/\eta}(K)\boxtimes_SK) \ar[r]^-{} & \delta^*_{X_\sigma/\sigma}\Psi_{X\times_SX}(\mathbb{D}_{X_\eta/\eta}(K)\boxtimes_SK) \ar[r]^-{} & \Psi_X\delta_{X_\eta/\eta}^*(\mathbb{D}_{X_\eta/\eta}(K)\boxtimes_SK) \\
  & \delta^*_{X_\sigma/\sigma}(\mathbb{D}_{X_\sigma/\sigma}(\Psi_X(K))\boxtimes_S\Psi_X(K))  \ar[ru]^-{} & \Psi_X(\mathbb{D}_{X_\eta/\eta}(K))\otimes\Psi_X(K). \ar[ru]^-{} \ar[l]_-{\sim} & 
    }
}
$$
\end{center}
While the other parts of the diagram above commute by functoriality, the commutativity of the upper left pentagon reduces to that of the following diagram:
\begin{align}
\begin{split}
  \xymatrix@=10pt{
  Th(-v_{|C_\sigma}) \ar[r]^-{} \ar[d]^-{} & i_C^*j_{C*}Th(-v_{|C_\eta}) \ar[d]^-{} \\
  \underline{Hom}(c_{\sigma1}^*\Psi_X(K),c_{\sigma1}^*\Psi_X(K)\otimes Th(-v_{|C_\sigma})) \ar[d]^-{} & \Psi_CTh(-v_{|C_\eta}) \ar[d]^-{} \\
  \underline{Hom}(c_{\sigma1}^*\Psi_X(K),\Psi_C(c_{\eta1}^*K)\otimes Th(-v_{|C_\sigma})) \ar[d]^-{u} & \Psi_C\underline{Hom}(c_{\eta1}^*K,c_{\eta1}^*K\otimes Th(-v_{|C_\eta})) \ar[d]^-{u} \ar[l]^-{} \\
  \underline{Hom}(c_{\sigma1}^*\Psi_X(K),\Psi_C(c_{\eta2}^!K)) \ar[d]^-{} & \Psi_C\underline{Hom}(c_{\eta1}^*K,c_{\eta2}^!K) \ar[ld]^-{} \\
  \underline{Hom}(c_{\sigma1}^*\Psi_X(K),c_{\sigma2}^!\Psi_C(K)) & 
  }
\end{split}
\end{align}
which follows from the naturality of the pseudo-monoidal structure of $\Psi$ in~\eqref{eq:pmonpsi}.
\endproof

\subsection{$\mathbb{A}^1$-homotopic analogue of a theorem of Verdier}
\label{sec:A1Verdier}
\subsubsection{}
The goal of this section is to prove the $\mathbb{A}^1$-homotopic analogue of a theorem of Verdier (\cite[Prop. 8.8]{Ver}), assuming that the exponential characteristic $p$ is invertible. The proof follows a pattern similar to \cite[\S3]{Var}, which is considerably simplified in our setting using $\mathbb{A}^1$-homotopies. 

\subsubsection{}
As the exponential characteristic $p$ is assumed to be invertible, by Lemma~\ref{lm:EK211}, we may replace the base field $k$ by its perfection and therefore suppose that $k$ is \emph{perfect}. Throughout this section we will make this assumption, which will not change the statements while making the proofs simpler.

\subsubsection{}
If $Z\to X$ is a closed immersion, let $D_ZX=Bl_{Z\times0}(X\times\mathbb{A}^1)-Bl_{Z\times0}(X\times0)$ be the (affine) deformation to the normal cone (\cite[\S5.1]{Ful}, \cite[10.3]{Ros}, \cite[3.2.3]{DJK}). Alternatively, $D_ZX$ can be defined as the spectrum over $\mathcal{O}_X$ of the Rees algebra
\begin{align}
\sum_n\mathcal{I}^n\cdot t^{-n}\subset\mathcal{O}_X[t,t^{-1}]
\end{align}
where $\mathcal{I}$ is the ideal sheaf defining $Z$ in $X$. 
There is a canonical morphism of schemes $D_ZX\to\mathbb{A}^1_k$, the fiber over $0$ is the normal cone $N_ZX=\operatorname{Spec}_{\mathcal{O}_X}(\oplus_{n\geqslant0}\mathcal{I}^n/\mathcal{I}^{n+1})$, while the fiber over $\mathbb{G}_{m}$ is simply $\mathbb{G}_{m,X}$. 
\begin{definition}
For any scheme $X$, denote by $q_X:\mathbb{G}_{m,X}\to X$ the projection. 
We define the \textbf{specialization to the normal cone} by 
\begin{align}
\label{eq:Verdiersp}
\begin{split}
sp_{X,Z}:\mathcal{SH}(X)&\to\mathcal{SH}(N_ZX)\\
K&\mapsto \Psi_{D_ZX}(q_X^*K).
\end{split}
\end{align}
\end{definition}

\subsubsection{}
\label{num:Var143}
(See also \cite[Lemma 1.4.3]{Var}) Let $f:Y\to X$ be a morphism of schemes, let $Z$ be a closed subscheme of $X$, and let $W$ be a closed subscheme of $f^{-1}(Z)$. Then the morphism $f$ lifts to morphisms $D_Z(f):D_WY\to D_ZX$ and $N_{W}(f):N_{W}{Y}\to N_ZX$. Indeed, let $\mathcal{I}'$ be the ideal sheaf defining $W$ in $Y$, then we have $f^\#(\mathcal{I})\subset \mathcal{I}'$ by assumption, which gives an inclusion $f^\#(\sum_n\mathcal{I}^n\cdot t^{-n})\subset \sum_n\mathcal{I}'^n\cdot t^{-n}$, and the result follows.

Consequently, the map~\eqref{eq:psiupper*} induces a natural transformation 
\begin{align}
\label{eq:BC^*}
(N_{W}(f))^*sp_{X,Z}\to sp_{Y,W}f^*. 
\end{align}
By adjunction, we have a natural transformation 
\begin{align}
\label{eq:BC_*}
sp_{X,Z}f_*\to (N_{W}(f))_*sp_{Y,W}.
\end{align}
If $f$ is smooth and $W=f^{-1}(Z)$, then the map~\eqref{eq:BC^*} is invertible. If $f$ is proper and $W=f^{-1}(Z)$, then the map~\eqref{eq:BC_*} is invertible. 

\subsubsection{}
\label{num:spethom}
It follows from~\ref{eq:spthom} that the specialization functor~\eqref{eq:Verdiersp} commutes with the formation of Thom spaces of vector bundles over $X$: if $v$ is a virtual vector bundle over $X$, then there is a canonical isomorphism
\begin{align}
sp_{X,Z}(K\otimes\operatorname{Th}(v))
\simeq 
sp_{X,Z}(K)\otimes\operatorname{Th}((v_{|D_ZX})_{|N_ZX}).
\end{align}
In particular, the specialization functor commutes with Tate twists.
\subsubsection{}
In the trivial case $Z=X$, we have a canonical isomorphism $D_ZX\simeq\mathbb{A}^1_Z$. The following lemma identifies the specialization map in this case:
\begin{lemma}
\label{lm:Vertriv}
For $Z=X$ and for any $K\in\mathcal{SH}(Z)[1/p]$, the canonical composition 
\begin{align}
\label{eq:Vertriv}
\begin{split}
K\to i^*j_*q_{Z}^*K\to \Psi_{\mathbb{A}^1_Z}(q_Z^*K)=sp_{Z,Z}(K)
\end{split}
\end{align}
is an isomorphism.

\end{lemma}

\proof
By Lemma~\ref{lm:CD7.2} and~\ref{num:spethom}, we may assume that $K=p_*\mathbbold{1}_Y$, where $p:Y\to Z$ is a proper morphism with $f:Y\to\operatorname{Spec}(k)$ smooth. Let $p':\mathbb{G}_m\times Y\to\mathbb{G}_m\times Z$ be the base change of $p$. 
In this case we have an isomorphism
\begin{align}
\label{eq:Vertrivnew}
\Psi_{\mathbb{A}^1_Z}(q_Z^*p_*\mathbbold{1}_Y)
\simeq 
\Psi_{\mathbb{A}^1_Z}(p'_*\mathbbold{1}_{\mathbb{G}_m\times Y})
\overset{\eqref{eq:psilower*}}{\simeq}
p_*\Psi_{\mathbb{A}^1_Y}\mathbbold{1}_{\mathbb{G}_m\times Y}
\overset{\eqref{eq:Ayo3510}}{\simeq}
p_*f^*\Psi_{\mathbb{A}^1_k}\mathbbold{1}_{\mathbb{G}_{m,k}}\simeq p_*\mathbbold{1}_Y.
\end{align}
By the naturality of $\Psi$, there is a canonical homotopy between the maps~\eqref{eq:Vertrivnew} and~\eqref{eq:Vertriv}, which proves the result.
\endproof

\subsubsection{}
For a closed immersion $Z\to X$ and $K\in\mathcal{SH}(X)[1/p]$, by~\ref{num:Var143} and Lemma~\ref{lm:Vertriv} we have a canonical map
\begin{align}
\label{eq:vermap}
sp_{X,Z}(K)_{|Z}\to sp_{Z,Z}(K_{|Z})\simeq K_{|Z}.
\end{align}
Our goal is to show Proposition~\ref{prop:Ver} below which states that the map~\eqref{eq:vermap} is always an isomorphism. Following the approach in \cite[\S3.1]{Var}, we first deal with the case of a strong deformation retract.
\begin{definition}
Let $Z\to X$ be a closed immersion. A \textbf{strong deformation} of $X$ onto $Z$ is a morphism $H:X\times\mathbb{A}^1\to X$ such that 
\begin{enumerate}
\item
$H_{|X\times\{1\}}$ is identity;
\item
$H_{|Z\times\mathbb{A}^1}:Z\times\mathbb{A}^1\to X$ agrees with the composition of the projection $\pi:Z\times\mathbb{A}^1\to Z$ and the inclusion $Z\to X$;
\item
$H_{|X\times\{0\}}$ is scheme-theoretically contained in $Z$. 
\end{enumerate}
In this case we say that $Z$ is a \textbf{strong deformation retract} of $X$.
\end{definition}

\subsubsection{}
For example, $\{0\}\to\mathbb{A}^1$ is a strong deformation retract via the morphism $H(x,a)=ax$. Similarly, the same result holds for $\mathbb{A}^{n-1}\times\{0\}\to\mathbb{A}^n$.

\begin{lemma}
\label{lm:var313}
If $Z$ is a strong deformation retract of $X$, then the map~\eqref{eq:vermap} is an isomorphism.
\end{lemma}

\proof
Let $H$ be a strong deformation of $X$ onto $Z$. Then $H$ induces a morphism $N_{\mathbb{A}^1_Z}H:N_{\mathbb{A}^1_Z}\mathbb{A}^1_X\to N_ZX$, with a commutative diagram
\begin{align}
\begin{split}
  \xymatrix@=11pt{
    \mathbb{A}^1_Z \ar[r]^-{} \ar[d]_-{} & Z \ar[d]^-{}\\
    N_{\mathbb{A}^1_Z}\mathbb{A}^1_X \ar[r]^-{N_{\mathbb{A}^1_Z}H} & N_ZX
  }
\end{split}
\end{align}
where the vertical maps are zero sections of the normal cone, and the upper horizontal map is the canonical projection. Then by~\ref{num:Var143}, the morphism $H$ induces a map
\begin{align}
\phi:
\pi^*sp_{X,Z}(K)_{|Z}
=
((N_{\mathbb{A}^1_Z}H)^*sp_{X,Z}(K))_{|\mathbb{A}^1_Z}
\to
(sp_{\mathbb{A}^1_X,\mathbb{A}^1_Z}(K_{|\mathbb{A}^1_X}))_{|\mathbb{A}^1_Z}
\simeq
\pi^*sp_{X,Z}(K)_{|Z}.
\end{align}
which is an endomorphism of the object $\pi^*sp_{X,Z}(K)_{|Z}\in\mathcal{SH}(\mathbb{A}^1_Z)$. 

For each point $a\in\mathbb{A}^1_k$, the composition of $H$ with the inclusion $X\simeq X\times\{a\}\to X\times\mathbb{A}^1$ gives a morphism $H_a:X\to X$, which induces a morphism $N_{Z}H_a:N_ZX\to N_ZX$. By~\eqref{eq:BC^*}, the morphism $H_a$ induces a map
\begin{align}
\phi_a:
sp_{X,Z}(K)_{|Z}
=
((N_{Z}H_a)^*sp_{X,Z}(K))_{|Z}
\to
sp_{X,Z}(K)_{|Z}.
\end{align}
The fiber of $\phi$ over $a$ is then $\phi_a$, so $\phi$ gives an $\mathbb{A}^1$-homotopy between $\phi_0$ and $\phi_1$. We know that $\phi_1$ is the identity map, while $\phi_0$ factors through the canonical map $sp_{X,Z}(K)_{|Z}\to K_{|Z}$ in~\eqref{eq:vermap}, and the result follows.
\endproof

\subsubsection{}
Now we are ready to prove the general case:
\begin{proposition}[Verdier's theorem on restriction to vertices]
\label{prop:Ver}
Let $K\in\mathcal{SH}(X)[1/p]$. Then the map $sp_{X,Z}(K)_{|Z}\to K_{|Z}$ in~\eqref{eq:vermap} is an isomorphism.
\end{proposition}
\proof
By Lemma~\ref{lm:CD7.2} and~\ref{num:spethom}, we may assume that $K=p_*\mathbbold{1}_Y$ where $p:Y\to X$ is a proper morphism with $Y$ smooth over $k$, such that the reduced scheme associated to $p^{-1}(Z)$ is a strict normal crossing divisor in $Y$. Let $q:p^{-1}(Z)\to Z$ be the morphism induced by $p$. Then we have $(p_*\mathbbold{1}_Y)_{|Z}\simeq q_*\mathbbold{1}_{p^{-1}(Z)}$, and $sp_{X,Z}(K)_{|Z}\simeq q_*sp_{Y,p^{-1}(Z)}(\mathbbold{1}_Y)$. Therefore we are reduced to prove the result for the pair $(Y,p^{-1}(Z))$ instead of the pair $(X,Z)$, so we may assume that $Z$ is a Cartier divisor in $X$.

Since the statement is local in the pair $(X,Z)$, by shrinking $X$ we may assume that $X$ is affine and $Z$ is defined by one equation in $X$, that is, there exists an integer $n$ and a closed immersion $i:X\to\mathbb{A}^n$ such that $Z=i^{-1}(\mathbb{A}^{n-1}\times\{0\})$. Denote by $i':Z\to\mathbb{A}^{n-1}\times\{0\}$ the closed immersion, then it suffices to show that the map
\begin{align}
\label{eq:newvermap}
i'_*(sp_{X,Z}(K)_{|Z})\to i'_*(K_{|Z})
\end{align}
obtained by applying $i'_*$ to the map~\eqref{eq:vermap}, is an isomorphism. But the map~\eqref{eq:newvermap} is canonically identified with the map
\begin{align}
\label{eq:Anvermap}
sp_{\mathbb{A}^n,\mathbb{A}^{n-1}\times\{0\}}(i_*K)_{|\mathbb{A}^{n-1}\times\{0\}}\to (i_*K)_{|\mathbb{A}^{n-1}\times\{0\}}
\end{align}
which is an isomorphism by the case of a strong deformation retract proved in Lemma~\ref{lm:var313}.
\endproof

\subsection{Deformation of correspondences}
\label{sec:defcorr}

\subsubsection{}
\label{num:defcorr}
Let $c:C\to X\times_kX$ be a correspondence and let $Z$ be a closed subscheme of $X$. By~\ref{num:Var143}, $c$ lifts to a correspondence $D_Z(c):D_{c^{-1}(Z\times_kZ)}C\to D_ZX\times_{\mathbb{A}^1_k}D_ZX$.

Over $\mathbb{G}_m$, this is $c_{\mathbb{G}_m}=c\times id_{\mathbb{G}_m}:C\times\mathbb{G}_m\to (X\times_{k}X)\times\mathbb{G}_m$.

Over $0$, this is $N_Z(c):N_{c^{-1}(Z\times_kZ)}C\to N_ZX\times_{k}N_ZX$ constructed in~\ref{num:Var143}.
\begin{lemma}
\label{lm:contrfix}
\begin{enumerate}

\item 
\label{lm:contrfix1}
The correspondence $c$ is contracting near $Z$ if and only if $c$ stabilizes $Z$ and the image of $N_Z(c)_1$ is set-theoretically supported at the zero section $Z\subset N_ZX$.

\item 
\label{num:contrtriv}
There is a canonical closed immersion $D_{c'^{-1}(Z)}Fix(c)\to Fix(D_Z(c))$. If $c$ is contracting near $Z$, $Fix(c)$ is connected, and 
$c'(Fix(c))\cap Z$ is non-empty, then there is a canonical isomorphism $(Fix(c)\times\mathbb{A}^1)_{red}\simeq (Fix(D_Z(c)))_{red}$, that is, $Fix(D_Z(c))$ is a topologically trivial family over $\mathbb{A}^1_k$.
\end{enumerate}
\end{lemma}
\proof
Let $\mathcal{I}'$ be the ideal sheaf defining $c^{-1}(Z\times_SZ)$ in $C$. Then the map $N_Z(c)_1:N_{c^{-1}(Z\times_SZ)}C\to N_ZX$ is given by the map
\begin{align}
\oplus_{n\geqslant0} ((\mathcal{I})^n/(\mathcal{I})^{n+1})
\to
\oplus_{n\geqslant0} ((\mathcal{I}')^n/(\mathcal{I}')^{n+1})
\end{align}
induced by $c_1:C\to X$, and the first claim follows. The second claim is proved in \cite[Thm. 2.1.3 (b)]{Var}, using \cite[Cor. 1.4.5]{Var} and apply Lemma~\ref{lm:icrediso}. 
\endproof

\subsubsection{}
\label{num:Kcondfin}
We now turn back to the proof of Proposition~\ref{prop:tr=0}. Let $K\in\mathcal{SH}_c(X)[1/p]$. We further consider the following condition:
\begin{itemize}
\item
The object $K_{|X\times\mathbb{G}_m}\in\mathcal{SH}_c(X\times\mathbb{G}_m)[1/p]=\mathcal{SH}_c((D_ZX)_\eta)[1/p]$ satisfies the conditions of~\ref{num:propnearby}.
\end{itemize}

\subsubsection{Proof of Proposition~\ref{prop:tr=0}}
\label{sec:tr=0}
By Lemma~\ref{lm:USLA}, $K$ is USLA over $k$, which implies that the object $K_{|X\times\mathbb{G}_m}\in\mathcal{SH}_c(X\times\mathbb{G}_m)[1/p]$ is USLA over $\mathbb{G}_{m,k}$. 
Consider the deformation construction in~\ref{num:defcorr}. By Lemma~\ref{lm:BCtrace} and Lemma~\ref{lm:sptrace}, we have the following commutative diagram
$$
\resizebox{\textwidth}{!}{
  \xymatrix{
    Map(c_{1}^*K,c_{2}^!K\otimes Th(v)) \ar[rr]^-{Tr(-/k)} \ar[d]_-{\eqref{eq:BCcorr}} & & H(Fix(c)/k,-v_{|Fix(c)}) \ar[d]^-{\eqref{eq:BCbiv}}\\
    Map((c_{\mathbb{G}_m})_{1}^*K_{|X\times\mathbb{G}_m},(c_{\mathbb{G}_m})_{2}^!K_{|X\times\mathbb{G}_m}\otimes Th(v_{|C\times\mathbb{G}_m})) \ar[rr]^-{Tr(-/\mathbb{G}_{m,k})} \ar[d]_-{\eqref{eq:spcorr}} & & H(Fix(c)\times\mathbb{G}_m/\mathbb{G}_{m,k},-v_{|Fix(c)\times\mathbb{G}_m}) \ar[d]^-{\eqref{eq:spFix}}\\
    Map(N_Z(c)_1^*sp_{X,Z}(K),N_Z(c)_2^!sp_{X,Z}(K)\otimes Th((v_{|D_Z(c)})_{|N_Z(c)})) \ar[rr]^-{Tr(-/k)} & & H(Fix(N_Z(c))/k,-v_{|Fix(N_Z(c))}).
  }
}
$$
By Lemma~\ref{lm:contrfix}~\eqref{num:contrtriv} and Lemma~\ref{lm:toptriv}, the composition of the two right vertical maps 
\begin{align}
H(Fix(c)/k,-v_{|Fix(c)})\to H(Fix(N_Z(c))/k,-v_{|Fix(N_Z(c))})
\end{align}
is an isomorphism. By Proposition~\ref{prop:Ver} we have $sp_{X,Z}(K)_{|Z}\simeq K_{|Z}\simeq0$, and by Lemma~\ref{lm:contrfix}~\eqref{lm:contrfix1} the image of $N_Z(c)_1$ is set-theoretically supported at the zero section $Z\subset N_ZX$, which implies that $N_Z(c)_1^*sp_{X,Z}(K)\simeq 0$. We conclude that the map
\begin{align}
Map(c_{1}^*K,c_{2}^!K\otimes Th(v))\xrightarrow{Tr(-/k)} H(Fix(c)/k,-v_{|Fix(c)})
\end{align}
is null-homotopic, which finishes the proof.
\endproof


\begin{thebibliography}{}

\bibitem[Ayo07a]{Ayo}
J. Ayoub,
\emph{Les six op\'erations de Grothendieck et le formalisme des cycles \'evanescents dans le monde motivique}, Ast\'erisque No. \textbf{314-315} (2007).

\bibitem[AIS17]{AIS}
J. Ayoub, F. Ivorra, J. Sebag, 
\emph{Motives of rigid analytic tubes and nearby motivic sheaves},  Ann. Sci. \'Ec. Norm. Sup\'er. (4) \textbf{50} (2017), no. 6, 1335-1382. 

\bibitem[Ayo14]{Ayo2}
J. Ayoub,
\emph{La r\'ealisation \'etale et les op\'erations de Grothendieck}, 
Ann. Sci. \'Ec. Norm. Sup\'er. (4) \textbf{47} (2014), no. 1, 1-145.

\bibitem[BH21]{BH}
T. Bachmann, M. Hoyois, 
\emph{Norms in motivic homotopy theory}, 
 Ast\'erisque No. 425 (2021).

\bibitem[BW21]{BW}
T. Bachmann, K. Wickelgren, 
\emph{$\mathbb{A}^1$-Euler classes: six functors formalisms, dualities, integrality and linear subspaces of complete intersections}, 
to appear in J. Inst. Math. Jussieu.

\bibitem[BD17]{BD}
M. Bondarko, F. D\'eglise,
\emph{Dimensional homotopy t-structures in motivic homotopy theory}, Adv. Math. \textbf{311} (2017), 91-189.

\bibitem[BBM+21]{BBMMO}
T. Brazelton, R. Burklund, S. McKean, M. Montoro, M. Opie,
\emph{The trace of the local $\mathbb{A}^1$-degree}, 
Homology Homotopy Appl. \textbf{23} (2021), no. 1, 243-255. 

\bibitem[BMP21]{BMP}
T. Brazelton, S. McKean, S. Pauli,
\emph{B\'ezoutians and the $\mathbb{A}^1$-degree}, 
\href{https://arxiv.org/abs/2103.16614}{arXiv:2103.16614}.

\bibitem[Cis13]{Cis1}
D.-C. Cisinski,
\emph{Descente par \'eclatements en K-th\'eorie invariante par homotopie},
Ann. of Math. (2) \textbf{177} (2013), no. 2, 425-448.

\bibitem[Cis21]{Cis}
D.-C. Cisinski,
\emph{Cohomological methods in intersection theory},
in \emph{Homotopy theory and arithmetic geometry – motivic and Diophantine aspects}, LMS-CMI research school, London, UK, July 9–13, 2018. Lecture notes. Cham: Springer. Lect. Notes Math. 2292, 49-105 (2021). 

\bibitem[CD15]{CD1}
D.-C. Cisinski, F. D\'eglise,
\emph{Integral mixed motives in equal characteristic}, 
Doc. Math. 2015, Extra vol.: Alexander S. Merkurjev's sixtieth birthday, 145-194.

\bibitem[CD16]{CD2}
D.-C. Cisinski, F. D\'eglise,
\emph{\'Etale motives}, 
Compos. Math. \textbf{152} (2016), no. 3, 556-666.

\bibitem[CD19]{CD}
D.-C. Cisinski, F. D\'eglise,
\emph{Triangulated categories of motives}, 
Springer Monographs in Mathematics. Springer, Cham, 2019.

\bibitem[DFJK21]{DFJK}
F. D\'eglise, J. Fasel, F. Jin, A. Khan, 
\emph{On the rational motivic homotopy category}, 
J. Ec. Polytech. Math. \textbf{8} (2021), 533-583.

\bibitem[DJK21]{DJK}
F. D\'eglise, F. Jin, A. Khan,
\emph{Fundamental classes in motivic homotopy theory},
J. Eur. Math. Soc. \textbf{23} (2021), no. 12, 3935-3993.

\bibitem[DG20]{DG}
B. Drew, M. Gallauer,
\emph{The Universal Six-Functor Formalism}, 
\href{https://arxiv.org/abs/2009.13610}{arXiv:2009.13610}.

\bibitem[EGA4$_{\textrm{III}}$]{EGA43}
A. Grothendieck, 
\emph{\'El\'ements de g\'eom\'etrie alg\'ebrique. III. \'Etude locale des sch\'emas et des morphismes de sch\'emas, Troisi\`eme partie}, 
r\'edig\'es avec la collaboration de Jean Dieudonn\'e,
Inst. Hautes \'Etudes Sci. Publ. Math. No. \textbf{28}, 1-255 (1966).

\bibitem[EGA4$_{\textrm{IV}}$]{EGA4}
A. Grothendieck, 
\emph{\'El\'ements de g\'eom\'etrie alg\'ebrique. IV. \'Etude locale des sch\'emas et des morphismes de sch\'emas, Quatri\`eme partie}, 
r\'edig\'es avec la collaboration de Jean Dieudonn\'e,
Inst. Hautes \'Etudes Sci. Publ. Math. No. \textbf{32}, 1-361 (1967).

\bibitem[EK20]{EK}
E. Elmanto, A. Khan,
\emph{Perfection in motivic homotopy theory},  Proc. Lond. Math. Soc. \textbf{120} (2020), no. 1, 28-38.

\bibitem[FHM03]{FHM}
H. Fausk, P. Hu, J. P. May,
\emph{Isomorphisms between left and right adjoints}, Theory Appl. Categ. \textbf{11} (2003), No. 4, 107-131.

\bibitem[Fer05]{Fer}
D. Ferrand,
\emph{On the non additivity of the trace in derived categories}, \href{https://arxiv.org/abs/math/0506589}{arXiv:math/0506589}.

\bibitem[Fuj97]{Fuj}
K. Fujiwara,
\emph{Rigid geometry, Lefschetz-Verdier trace formula and Deligne’s conjecture},
Invent. Math. \textbf{127}, No. 3, 489-533 (1997). 

\bibitem[Ful98]{Ful}
W. Fulton,
\emph{Intersection theory}, Second edition, Ergebnisse der Mathematik und ihrer Grenzgebiete. 3. Folge. A Series of Modern Surveys in Mathematics, \textbf{2}. Springer-Verlag, Berlin, 1998.

\bibitem[GM93]{GM}
M. Goresky, R. MacPherson,
\emph{Local contribution to the Lefschetz fixed point formula},
Invent. Math. \textbf{111}, No. 1, 1-33 (1993). 

\bibitem[Gro66]{Gro}
A. Grothendieck,
\emph{Formule de Lefschetz et rationalit\'e des fonctions $L$},
S\'emin. Bourbaki Vol. 9, 17e année (1964/1965), Exp. No. 279, 15 p. (1966). 

\bibitem[Hoy15]{Hoy}
M. Hoyois,
\emph{A quadratic refinement of the Grothendieck-Lefschetz-Verdier trace formula}, Algebraic \& Geometric Topology \textbf{14} (2015), no.~6, 3603-3658.

\bibitem[Ill17]{Ill}
L. Illusie,
\emph{Around the Thom-Sebastiani theorem, with an appendix by Weizhe Zheng}, Manuscripta Math. \textbf{152} (2017), no. 1-2, 61-125.

\bibitem[ILO14]{ILO}
L. Illusie, Y. Laszlo, F. Orgogozo,
\emph{Travaux de Gabber sur l'uniformisation locale et la cohomologie \'etale des sch\'emas quasi-excellents}, 
S\'eminaire \`a l'Ecole Polytechnique 2006-2008. With the collaboration of Fr\'ed\'eric D\'eglise, Alban Moreau, Vincent Pilloni, Michel Raynaud, Jo\"el Riou, Beno\^it Stroh, Michael Temkin and Weizhe Zheng. Ast\'erisque No. 363-364, Soci\'et\'e Math\'ematique de France, Paris, 2014.

\bibitem[JY21a]{JY}
F. Jin, E. Yang,
\emph{K\"unneth formulas for motives and additivity of traces}, 
Adv. Math. \textbf{376} (2021), Article ID 107446.

\bibitem[JY21b]{JY2}
F. Jin, E. Yang,
\emph{Some results on the motivic nearby cycle}, \href{https://arxiv.org/abs/2107.08603}{arXiv:2107.08603}.

\bibitem[KW19]{KW}
J. Kass, K. Wickelgren,
\emph{The class of Eisenbud-Khimshiashvili-Levine is the local $\mathbb{A}^1$-Brouwer degree}, Duke Math. J. \textbf{168} (2019), no. 3, 429-469. 

\bibitem[Kha16]{Kha}
A. Khan,
\emph{Motivic homotopy theory in derived algebraic geometry}, 
Ph.D. thesis, Universit\"at Duisburg-Essen, 2016, available at https://www.preschema.com/thesis/.

\bibitem[Lev20]{Lev}
M. Levine,
\emph{Aspects of enumerative geometry with quadratic forms}, 
Doc. Math. \textbf{25} (2020), 2179-2239.

\bibitem[LR20]{LR}
M. Levine, A. Raksit,
\emph{Motivic Gauss-Bonnet formulas}, 
Algebra Number Theory, \textbf{14} (7):1801-1851, 2020.

\bibitem[LYZR19]{LYZR}
M. Levine, Y. Yang, G. Zhao, J. Riou,
\emph{Algebraic elliptic cohomology theory and flops I}, 
Math. Ann. \textbf{375} (2019), no. 3-4, 1823-1855.

\bibitem[LZ22]{LZ}
Q. Lu, W. Zheng,
\emph{Categorical traces and a relative Lefschetz–Verdier formula},
Forum of Mathematics, Sigma, \textbf{10}, E10, 2022.

\bibitem[Mor12]{Mor}
F. Morel,
\emph{$\mathbb{A}^1$-algebraic topology over a field}, 
Lecture Notes in Mathematics, 2052. Springer, Heidelberg, 2012.

\bibitem[Ols16]{Ols}
M. Olsson,
\emph{Motivic cohomology, localized Chern classes, and local terms}, Manuscripta Math. \textbf{149} (2016), no. 1-2, 1-43.

\bibitem[Pin92]{Pin}
R. Pink,
\emph{On the calculation of local terms in the Lefschetz-Verdier trace formula and its application to a conjecture of Deligne},
Ann. Math. (2) \textbf{135}, No. 3, 483-525 (1992). 

\bibitem[Ros96]{Ros}
M. Rost,
\emph{Chow groups with coefficients},
Doc. Math. \textbf{1}, 319-393 (1996). 

\bibitem[SGA5]{SGA5}
A. Grothendieck,
\emph{Cohomologie l-adique et fonctions L}, S{\'e}minaire de g{\'e}om{\'e}trie alg{\'e}brique du Bois-Marie 1965-66 (SGA 5). Avec la collaboration de I. Bucur, C. Houzel, L. Illusie, J.-P. Jouanolou, et J.-P. Serre. Springer Lecture Notes, Vol. 589. Springer-Verlag, Berlin-New York, 1977.

\bibitem[SGA6]{SGA6}
P. Berthelot, A. Grothendieck, L. Illusie,
\emph{Th\'eorie des intersections et th\'eor\`eme de Riemann-Roch}, S\'eminaire de G\'eom\'etrie Alg\'ebrique du Bois-Marie 1966-1967 (SGA 6).
With the collaboration of D. Ferrand, J. P. Jouanolou, O. Jussila, S. Kleiman, M. Raynaud et J. P. Serre.  Lecture Notes in Mathematics, Vol. \textbf{225}. Springer-Verlag, Berlin-New York, 1971.

\bibitem[Var07]{Var}
Y. Varshavsky,
\emph{Lefschetz-Verdier trace formula and a generalization of a theorem of Fujiwara},
Geom. Funct. Anal. \textbf{17}, No. 1, 271-319 (2007). 

\bibitem[Ver81]{Ver}
J.-L. Verdier,
\emph{Sp\'ecialisation de faisceaux et monodromie mod\'er\'ee},
in \emph{Analysis and topology on singular spaces} (Luminy, 1981), 332-364,
Ast\'erisque, \textbf{101-102}, Soc. Math. France, Paris, 1983.

\bibitem[YZ21]{YZ}
E. Yang, Y. Zhao,
\emph{On the Relative Twist Formula of $\ell$-adic Sheaves},
Acta Math. Sin. (Engl. Ser.) \textbf{37} (2021), no. 1, 73-94. 

\end{thebibliography}
\end{document}